\documentclass[a4paper,10pt]{article}

\usepackage{amsthm,amssymb,amsfonts,amscd,amsmath,mathtools}
\usepackage{stmaryrd}
\usepackage{csquotes}
\usepackage{hyperref}
\usepackage{microtype}
\usepackage{tikz}
\usepackage{xcolor}
\usepackage{soul}

\usetikzlibrary{decorations.pathreplacing,matrix,calc,arrows,shapes.geometric,arrows.meta}

\definecolor{myblue}{rgb}{0,0,0.6}
\hypersetup{pdftitle={Trefftz Discontinuous Galerkin Approximation of the T-Matrix for Scattering by Periodic Layered Structures},  pdfauthor={Moiola, Monforte, Zanotto},    pdffitwindow=true,
   colorlinks=true, linkcolor=myblue,  citecolor=myblue, filecolor=myblue,   urlcolor=myblue,
   bookmarksopen=true}

\usepackage[backend=biber,isbn=false,giveninits=true,doi=false,sortcites,url=false,eprint=false]{biblatex}
\newcommand{\dlgraffa}{\{ \hspace{-0.1cm} \{}
\newcommand{\drgraffa}{\} \hspace{-0.1cm} \}}
\newcommand{\vertiii}[1]{{\left\vert\kern-0.25ex\left\vert\kern-0.25ex\left\vert #1 \right\vert\kern-0.25ex\right\vert\kern-0.25ex\right\vert}}

\newcommand{\IR}{{\mathbb R}}
\newcommand{\IC}{{\mathbb C}}
\newcommand{\IN}{{\mathbb N}}
\newcommand{\IZ}{{\mathbb Z}}
\newcommand{\bx}{{\mathbf x}}
\newcommand{\bd}{{\mathbf d}}
\newcommand{\bn}{{\mathbf n}}
\newcommand{\bE}{{\mathbf E}}
\newcommand{\bT}{{\mathbf T}}
\newcommand{\Mh}{{\mathcal M_h}}
\newcommand{\ri}{{\mathrm i}}
\newcommand{\re}{{\mathrm e}}
\newcommand{\rd}{{\,\mathrm d}}
\newcommand{\inc}{{\mathrm{inc}}}

\newcommand{\scat}{{\mathrm{scat}}}
\newcommand{\cA}{{\mathcal A}}
\newcommand{\pa}{{\mathtt{a}}}
\newcommand{\pb}{{\mathtt{b}}}
\newcommand{\pd}{{\mathtt{d}}}
\newcommand{\dn}{{\partial_\bn}}
\newcommand{\dxo}{{\partial_{x_1}}}
\newcommand{\dxt}{{\partial_{x_2}}}
\newcommand{\Gjp}{{\Gamma_j^+}}
\newcommand{\Gjm}{{\Gamma_j^-}}
\newcommand{\Gjpm}{{\Gamma_j^\pm}}
\newcommand{\Strip}{{\Upsilon}}
\newcommand{\aup}{\mathbf{a}^\uparrow}
\newcommand{\adown}{\mathbf{a}^\downarrow}
\newcommand{\bup}{\mathbf{b}^\uparrow}
\newcommand{\bdown}{\mathbf{b}^\downarrow}
\newcommand{\Rup}{\mathbf{R}^\uparrow}
\newcommand{\Rdown}{\mathbf{R}^\downarrow}
\newcommand{\Tup}{\mathbf{T}^\uparrow}
\newcommand{\Tdown}{\mathbf{T}^\downarrow}
\newcommand{\DtN}{{\mathrm{DtN}}}
\newcommand{\Phiup}[1]{{\Phi_{#1}^{\uparrow}}}
\newcommand{\Psidown}[1]{{\Psi_{#1}^{\downarrow}}}
\newcommand{\Psiup}[1]{{\Psi_{#1}^{\uparrow}}}
\newcommand{\Phidown}[1]{{\Phi_{#1}^{\downarrow}}}
\newcommand{\Tmat}{{\mathbb T}}
\newcommand{\bD}{{\mathbf D}}
\newcommand{\bI}{{\mathbf I}}
\newcommand{\bzero}{{\mathbf 0}}
\newcommand{\Gleft}{\Gamma_{\mathrm{left}}}
\newcommand{\Gright}{\Gamma_{\mathrm{right}}}
\newcommand{\GD}{{\Gamma_D}}
\newcommand{\dnK}{{\partial_{\bn_K}}}

\newcommand{\cv}{{\overline{v}}}
\newcommand{\cw}{{\overline{w}}}
\newcommand{\Fh}{\mathcal{F}_h}
\newcommand{\FhI}{\mathcal{F}_h^I}
\newcommand{\GHp}{{\Gamma_{H^+}}}
\newcommand{\GHm}{{\Gamma_{H^-}}}
\newcommand{\GHpm}{{\Gamma_{H^\pm}}}
\newcommand{\BoxMat}[4]{{\fbox{$\begin{matrix} #1 &\hspace*{-5pt} #2 \\ #3 &\hspace*{-5pt} #4 \end{matrix}$}}}

\graphicspath{{figs/}}

\title{Trefftz DG Approximation of the T-Matrix for Scattering by Periodic Layered Structures}
\author{Armando Maria Monforte\thanks{(Corresponding author) Department of Mathematics, University of Pavia, Via Ferrata 5, Pavia, Italy
(\href{mailto:armandomaria.monforte01@universitadipavia.it}{armandomaria.monforte01 @universitadipavia.it}), ORCID: 0009-0000-7687-2217}, 
Andrea Moiola\thanks{Department of Mathematics, University of Pavia, Via Ferrata 5, Pavia, Italy 
(\href{mailto:andrea.moiola@unipv.it}{andrea.moiola@unipv.it}), ORCID: 0000-0002-6251-4440},
Simone Zanotto\thanks{Istituto Nanoscienze – CNR, NEST-SNS, Piazza San Silvestro 12, Pisa, Italy
(\href{mailto:simone.zanotto@nano.cnr.it}{simone.zanotto@nano.cnr.it}), ORCID: 0000-0001-7180-3335}}
\date{\today}

\newtheorem{thm}{Theorem}[section]
\newtheorem{defin}[thm]{Definition}
\newtheorem{rem}[thm]{Remark}

\allowdisplaybreaks

\begin{document}
\maketitle

\section*{Abstract}
We study the scattering of time-harmonic electromagnetic waves by periodic layered gratings, modelled by the 2D Helmholtz equation. 
The periodic obstacle may include penetrable and impenetrable regions, and consists of a finite number of stacked layers.
The boundary value problem is formulated on a single periodic cell using quasi-periodic boundary conditions.
The radiation condition in the vertical directions is imposed through Dirichlet-to-Neumann (DtN) operators. 

To efficiently treat multilayer configurations, we adopt a formulation based on the T-matrix method.
The global scattering problem is decomposed into boundary value problems posed on individual layers.
On the layer boundaries, the field is expressed in terms of quasi-periodic modal expansions, and the layer T-matrix describes the map between incoming and outgoing wave modes.
Each local T-matrix is approximated numerically using a plane-wave based Trefftz Discontinuous Galerkin (TDG) method, which provides an efficient discretization of the layer scattering response.
The T-matrix technique leads to linear computational complexity in the number of layers in the grating.

\bigskip
\textbf{Keywords:} \ Diffraction grating, Quasi-periodic, Helmholtz equation, T-matrix, Discontinuous Galerkin, Trefftz method, Plane wave basis

\bigskip
\textbf{Mathematics Subject Classification (2020):} \ 
65N30, 35J05, 35Q60, 78A45, 78M10

\section{Introduction}

The scattering of time-harmonic electromagnetic waves by periodic structures, known as diffraction gratings, arises in a wide range of applications in photonics, optics, and metamaterial design
\cite{BaoLi2022,Joannopoulos2011}.
In the two-dimensional transverse-electric setting, Maxwell's equations reduce to the Helmholtz equation with piecewise-constant, periodic relative permittivity $\varepsilon$.
Classical numerical approximations of such diffraction gratings reduce the computational domain to a single periodic cell using quasi-periodic boundary conditions and Dirichlet-to-Neumann operators \cite{BonnetBD1994,Kirsch1993,Civiletti2020}.

To curb the computational cost of discretizing large and complex cells, we borrow the T-matrix technique from multiple-scattering theory \cite{Martin2006,Waterman1965,Mishchenko1996}.
To approximate the scattering by many obstacles, each of these is associated to a so-called T-matrix, which maps the Fourier coefficients of any wave impinging on the obstacle into the Fourier coefficients of the corresponding scattered field.
The size of a T-matrix is chosen in dependence of 
the wavenumber and the obstacle size, but independently 
of the obstacle shape and material properties, although a larger size may be needed for accuracy in presence of resonant structures exciting strongly evanescent modes.
Multiple scattering problems are
solved by constructing a linear system that couples the T-matrices of each obstacle, leading to considerable savings in the computational effort.

We adapt this idea to the quasi-periodic setting.
First, the domain $\Omega$ is split in $N$ layers $\Omega_j$, which may contain inhomogeneities and impenetrable regions.
Second, a T-matrix $\Tmat_j$, mapping incoming to outgoing fields, is computed for each layer, exploiting the natural Fourier expansion dictated by quasi-periodicity.
Third, the T-matrices of these layers are coupled through propagation matrices into a global banded block linear system.

A key feature of this formulation is the decoupling of scales, or model order reduction:
the global system is $N\times N$ block tridiagonal, with $(4M+2)\times(4M+2)$-sized blocks, 
where $M$ is the number of Fourier modes used to represent the outgoing and incoming waves on a layer, and $N$ is the number of layers (see \eqref{eq:coupling_system}).
These parameters are independent of the vertical size of the layers, their material properties, the shape of the obstacles contained, and also on the discretization and the number of degrees of freedom used inside each layer to approximate the T-matrix.
The local problems can therefore be discretized with a fine discretization or a high-order method without increasing the cost of the global coupling step.

This is not the first use of Fourier-basis model order reduction for the Helmholtz equation with scatterers partitioned in slabs.
A related method was proposed in \cite{MontielSquireBennetts2015} for circular sound-hard scatterers collected in unbounded slabs (without quasi-periodicity),
and a similar approach based on Neumann-to-Dirichlet matrices was proposed in \cite{Liu2026} for non-periodic metasurfaces divided into unit cells.

In this work we approximate each T-matrix using the Trefftz Discontinuous Galerkin (TDG) method with plane-wave basis functions \cite{Hiptmair2016,Kapita2018}, combined with the DtN boundary conditions on the horizontal layer interfaces.
We extend the method of \cite{Monforte2025} to allow for arbitrary quasi-periodic incident fields, which is necessary to build the T-matrix.

The proposed framework offers several advantages over a monolithic DtN-TDG solver on the full quasi-periodic domain: reusability of repeated layers, parallel assembly of the T-matrices, and scalability in the number of layers. Indeed, adding more layers to the structure does not strongly affect the computational cost, while
the direct DtN-TDG on the full domain requires the solution of a much larger linear system whose size is proportional to $N$.

In \S\ref{s:periodic_general} we briefly introduce the Helmholtz quasi-periodic boundary value problems (BVP) of interest, together with the associated DtN operators.
In \S\ref{s:LayeredStructures} we introduce the notation and the assumptions related to the layers making up the scatterer, and write the local BVPs that will be approximated numerically.
In \S\ref{s:Tmat} we explain in detail the application of the T-matrix idea to the quasi-periodic layered setting.
In \S\ref{s:TDG} we describe the DtN-TDG method used to approximate the T-matrices on each layer.
Finally, in \S\ref{s:Numerics} we show some numerical results obtained with an implementation of the proposed method, which is available online, including a comparison with the RCWA method on a physics-relevant example. There, we also show that the DtN-TDG method behaves well in the vicinity of a bound state in the continuum (BIC), i.e.~when the quality factor of a guided mode resonance (GMR) approaches infinity and the quasi-periodic BVP loses unique solvability.
In other words, the proposed method is able to approximately detect the presence of a Rayleigh--Bloch wave within the radiation continuum (i.e.~above the light line) \cite{Bennetts_Peter_2022, chaplain_acoustic_2025}.

\section{Periodic scattering problem} \label{s:periodic_general}
We present the problem of the scattering of a plane wave by a periodic structure, also called diffraction grating,
which has been the subject of significant research in the literature using different approaches, including \cite{BonnetBD1994,Civiletti2020,Kirsch1993} and the references therein.
We assume a time-harmonic dependence
of the form $\re^{-\ri\omega t}$, where $\ri = \sqrt{-1}$ and $\omega$ is the wave angular frequency. 

We consider a TE-polarized incident field, which leads to an electric field of the form
$\mathbf E(x_1,x_2,x_3)=(0,0,E_3(x_1,x_2))$, solving time-harmonic Maxwell equations in materials that are invariant in the direction $x_3$. Maxwell equations then reduce to the two-dimensional Helmholtz equation in the variables $x_1,x_2$ for the component $u=E_3$, \cite[eq.~(2.5)]{BonnetBD1994}.
TM-polarized waves can be treated similarly; see Remark~\ref{rem:TMmodes}.

\subsection{Domain description}\label{s:Domain}
We consider gratings that are periodic in the direction $x_1$ with period $L>0$. We allow both perfect electric conductor (PEC) obstacles, which are impenetrable by waves, and dielectric media. The PEC obstacle is represented by an open Lipschitz set $D\subset\IR^2$, possibly empty, that is periodic in the $x_1$-direction with period $L$. 
The scattering region is $\Omega_0:=\IR^2\setminus \overline D$, and we assume it to be connected (although the problem and the method proposed can easily be extended to the case with a PEC interface $\{\bx\in\IR^2:\,x_2=f_{\mathrm{PEC}}(x_1)\}$ for a Lipschitz $L$-periodic $f_{\mathrm{PEC}}$).
Let $\varepsilon_0,\mu_0>0$ denote the permittivity and permeability of vacuum. 
We assume that the relative magnetic permeability is constant $\mu=1$ in $\Omega_0$.
The relative permittivity of the medium $\varepsilon\in L^\infty(\Omega_0)$ is a complex-valued, piecewise-constant function with $\Re(\varepsilon) > 0$ and $\Im(\varepsilon) \geq 0$, which is periodic with period $L$ in $x_1$.
We assume that the region containing the inhomogeneities
of $\varepsilon$ is bounded in the $x_2$-direction, i.e.\  $\varepsilon(\bx)=\varepsilon^+$ 
in $\{x_2>H^+\}$ and $\varepsilon(\bx)=\varepsilon^-$ in $\{x_2<H^-\}$ for some $H^-<H^+\in\IR$.
We also assume that $\varepsilon^+$ is real and positive.

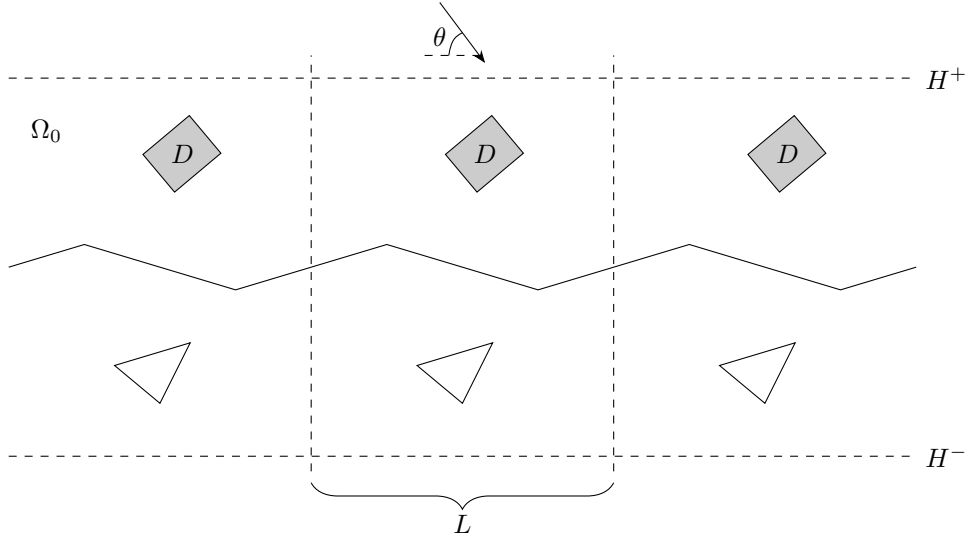
\begin{figure}[htbp]
\centering
\begin{tikzpicture}
\draw[dashed] (0,-0.2) -- (0,5.3);
\draw[dashed] (4,-0.2) -- (4,5.3);
\draw[dashed] (-4,5) -- (8,5);
\draw (8,5) node[right]{$H^+$};
\fill[gray!40!white, draw=black, rotate around={40:(-1.3,3.2)}] (-1.5,3.75) rectangle (-0.7,4.4);
\fill[gray!40!white, draw=black, rotate around={40:(2.7,3.2)}] (2.5,3.75) rectangle (3.3,4.4);
\fill[gray!40!white, draw=black, rotate around={40:(6.7,3.2)}] (6.5,3.75) rectangle (7.3,4.4);
\node[] at (-1.7,4) {$D$};
\node[] at (2.3,4) {$D$};
\node[] at (6.3,4) {$D$};
\draw (-4,2.5) -- (-3,2.8) -- (-1,2.2) -- (0,2.5) -- (1,2.8) -- (3,2.2) -- (4,2.5) -- (5,2.8) -- (7,2.2) -- (8,2.5);
\coordinate (A) at (1.4,1.2);
\coordinate (B) at (2,0.7);
\coordinate (C) at (2.4,1.5);
\draw (A) -- (B) -- (C) -- cycle;
\coordinate (A1) at (-2.6,1.2);
\coordinate (B1) at (-2,0.7);
\coordinate (C1) at (-1.6,1.5);
\draw (A1) -- (B1) -- (C1) -- cycle;
\coordinate (A2) at (5.4,1.2);
\coordinate (B2) at (6,0.7);
\coordinate (C2) at (6.4,1.5);
\draw (A2) -- (B2) -- (C2) -- cycle;
\draw[dashed] (-4,0) -- (8,0);
\draw (8,0) node[right]{$H^-$};
\draw[decorate, decoration={brace,amplitude=10}] (4,-0.35)--(0,-0.35);
\draw(2,-0.9)node{$L$};
\draw (1.82,5.3) .. controls (1.85,5.5) and (1.9,5.55) .. (2.0, 5.6);
\node[] at (1.7,5.55) {$\theta$};
\draw[-{Stealth[slant=0]}] (1.7,6) -- (2.3,5.2);
\draw[dashed] (1.5,5.3) -- (2.25,5.3);
\node[] at (-3.5,4.3) {$\Omega_0$};
\node[] at (1,4.5) {$\varepsilon^+$};
\node[] at (3,0.5) {$\varepsilon^-$};
\node[] at (1.9,1.1) {$\varepsilon_{in}$};
\node[] at (5.9,1.1) {$\varepsilon_{in}$};
\node[] at (-2.1,1.1) {$\varepsilon_{in}$};
\end{tikzpicture}  
\caption{Geometry of the periodic scattering region $\Omega_0$ and the Dirichlet obstacle $D$.
The permittivity $\varepsilon$ takes constant values $\varepsilon^+$, $\varepsilon^-$ and $\varepsilon_{in} \neq \varepsilon^-$ in the regions of $\Omega_0=\IR^2\setminus\overline D$ bounded by the continuous lines.
}
\label{fig:scatt_layered}
\end{figure}

A possible scattering region is depicted in Figure \ref{fig:scatt_layered}: shaded regions represents the impenetrable obstacle $D$, while continuous lines represent the interfaces separating different materials, i.e.\ regions with constant $\varepsilon$.

We define a truncated domain taking into account the periodicity in $x_1$ and the boundedness in $x_2$ of the scatterer:
\begin{equation}\label{eq:Omega}
\Omega:=\{\bx\in\Omega_0:\; 0<x_1<L,\; H^-<x_2<H^+\}.
\end{equation}
We introduce also some notation for the parts of the boundary of $\Omega$ and the upper/lower regions:
\begin{align*}
\GHpm:=&\{\bx\in\IR^2: 0\le x_1\le L,\; x_2= H^\pm \},\\
\Gleft:=&\{\bx\in\IR^2:  x_1=0,\; H^- \le x_2\le H^+ \},\\
\Gright:=&\{\bx\in\IR^2:  x_1=L,\; H^- \le x_2\le H^+ \},\\
\GD:=& \partial D\cap\overline\Omega,\\
\Strip^+:=&\{\bx\in\IR^2:   0<x_1<L,\; x_2> H^+ \},\\
\Strip^-:=&\{\bx\in\IR^2:   0<x_1<L,\; x_2< H^- \},
\end{align*}
so that $\partial\Omega=\GHp\cup\GHm\cup\Gleft\cup\Gright\cup\GD$ and $\varepsilon|_{\Strip^\pm}=\varepsilon^\pm$.
We denote by $\bn$ the outward-pointing unit normal on $\partial\Omega$.

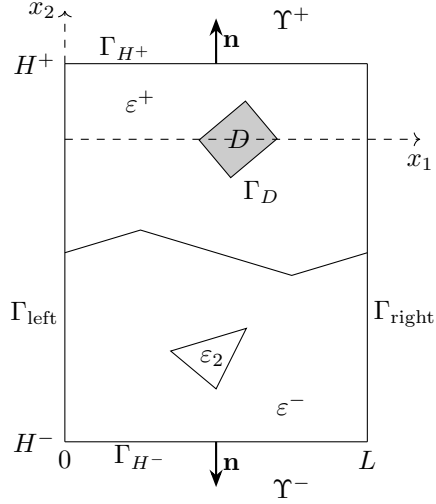
\begin{figure}[htbp]
\centering
\begin{tikzpicture}
\draw[] (0,5) -- (4,5);
\draw[] (0,0) -- (4,0);
\draw[] (0,0) -- (0,5);
\draw[] (4,0) -- (4,5);
\draw (0,5) node[left]{$H^+$};
\draw (0,0) node[left]{$H^-$};
\draw (0,0) node[below]{$0$};
\draw (4,0) node[below]{$L$};
\fill[gray!40!white, draw=black, rotate around={40:(2.7,3.2)}] (2.5,3.75) rectangle (3.3,4.4);
\node[] at (2.3,4) {$D$};
\draw (0,2.5) -- (1,2.8) -- (3,2.2) -- (4,2.5);
\coordinate (A) at (1.4,1.2);
\coordinate (B) at (2,0.7);
\coordinate (C) at (2.4,1.5);
\draw (A) -- (B) -- (C) -- cycle;
\node[] at (4.7,3.7) {$x_1$};
\node[] at (-0.3,5.7) {$x_2$};
\node[] at (-0.4,1.7) {$\Gleft$};
\node[] at (4.5,1.7) {$\Gright$};
\node[] at (0.8, 5.2) {$\GHp$};
\node[] at (1, -0.2) {$\GHm$};
\node[] at (1,4.5) {$\varepsilon^+$};
\node[] at (2.6,3.3) {$\Gamma_D$};
\node[] at (3,0.5) {$\varepsilon^-$};
\draw [dashed,->](0,5) -- (0,5.7); 
\draw [dashed,->](0,4) -- (4.7,4); 
\draw [-{Stealth[slant=0]}, thick](2,5) -- (2,5.6);
\node[] at (2.2,5.3) {$\mathbf{n}$};
\draw [-{Stealth[slant=0]}, thick](2,0) -- (2,-0.6);
\node[] at (2.2,-0.3) {$\mathbf{n}$};
\node[] at (3,5.6) {$\Strip^+$};
\node[] at (3,-.6) {$\Strip^-$};
\node[] at (1.9,1.1) {$\varepsilon_{in}$};
\end{tikzpicture}  
\caption{The truncated domain $\Omega=(0,L)\times(H^-,H^+)\setminus \overline D$. 
The relative permittivity takes the values $\varepsilon^+$, $\varepsilon^-$ and $\varepsilon_{in} \neq \varepsilon^-$ in the three regions separated by continuous lines.}
\label{fig:omega_layered}
\end{figure}

In Figure \ref{fig:omega_layered} we see the domain $\Omega$ obtained from the truncation of the region $\Omega_0$ in Figure~\ref{fig:scatt_layered}.

\subsection{Helmholtz equation and quasi-periodic conditions}
The source of the scattering problem is a 
plane wave:
\begin{equation}\label{eq:Incident}
\bE^\inc:=(0,0,u^\inc),\qquad
u^\inc(\bx):=u^\inc(x_1,x_2)=\exp\{\ri\kappa^+(x_1 \cos \theta + x_2 \sin \theta)\},
\end{equation} 
where $\theta \in (-\pi, 0)$ is the wave propagation angle against the horizontal, $\kappa^+:= k\sqrt{\varepsilon^+}$ is the wavenumber in the upper region $\Strip^+$ and $k$ is the wavenumber in free space.
We call the plane wave with $\theta \in (-\pi, 0)$ ``downward-propagating'' as the second component of the propagation direction is negative. 
The incident plane wave direction is in the $\{x_1,x_2\}$-plane and hence normal to the periodic structure.
We define the piecewise-constant wavenumber function $\kappa(\bx):=k\sqrt{\varepsilon(\bx)}\in L^\infty(\Omega_0)$.

We denote by $u$ the third component of the total electric field $\mathbf E=(0,0,u)$ generated by the scattering of $\bE^\inc$ on $D$ and on the heterogeneities of $\varepsilon$.
From the Maxwell equations, it is classical (e.g.\ \cite{BonnetBD1994,Civiletti2020}) that $u$ satisfies the Helmholtz equation 
\begin{equation*}
\Delta u+k^2\varepsilon u=0 \qquad \text{in }\Omega_0,
\end{equation*}
which has to be complemented with appropriate boundary, radiation and periodicity condition.

We impose the boundary condition $u=0$ on $\partial D$, since we are assuming that $D$ is a perfect electric conductor (PEC).

The incoming wave $u^\inc$ is periodic in the direction $x_1$ with period $2\pi/(\kappa^+\cos\theta)$, which is in general different from $L$.
However, we observe that $u^\inc$ satisfies the relation $u^\inc(x_1+L,x_2)=\re^{\ri\kappa^+L\cos\theta}u^\inc(x_1,x_2)$.
This suggests the following classical definition \cite{Pinto2020}.

\begin{defin}[Quasi-periodic function]\label{def:QP}
A function $u\in \mathcal{C}^0(\IR^2)$ is called quasi-periodic in $\IR^2$, of period $L$, with parameter $\alpha_0\in\IR$, if
\begin{equation*}
u(x_1+L,x_2)=\re^{\ri\alpha_0 L}u(x_1,x_2) \qquad \forall \bx = (x_1,x_2) \in \IR^2.
\end{equation*}
A function $u\in \mathcal{C}^0(\Omega_0)$ is quasi-periodic in $\Omega_0$ if the same relation holds for all $\bx\in\Omega_0$.
\end{defin} 

The total field $u$ inherits the quasi-periodicity with parameter $\alpha_0=\kappa^+\cos\theta$ from the plane wave $u^\inc$; note that $\alpha_0$ is also referred to as the Bloch wavenumber (or Floquet wavenumber). 
To enforce this property, we restrict the Helmholtz equation to the bounded domain $\Omega$ in \eqref{eq:Omega} and impose a quasi-periodicity condition on $\Gleft$ and $\Gright$:
\begin{align*}
\begin{aligned}
u(L,x_2) =&\ \re^{\ri\alpha_0 L}u(0,x_2) \\ 
\dxo  u(L,x_2) =&\ \re^{\ri\alpha_0 L}\dxo  u (0,x_2) 
\end{aligned}
\qquad H^- < x_2 < H^+, \; (0,x_2)\notin\overline D.
\end{align*}

On $\GHp$ and $\GHm$ we impose a radiation condition, asking that the scattered field propagates upward on $\GHp$ and the total field propagates downward on $\GHm$; this is equivalent to asking that these two boundaries are transparent to waves propagating away from the grating. A way to impose the radiation condition is using Dirichlet-to-Neumann operators, which we introduce in the following section.

\subsection{Dirichlet-to-Neumann operators} \label{s:DtN}
We first introduce the quasi-periodic Sobolev spaces on $\Omega$ and on the horizontal boundaries $\GHpm$, which are used to define the DtN operators.
Following \cite[\S3.1]{BonnetBD1994}, for $\alpha_0\in\IR$, the space $\mathcal{C}_{\alpha_0}^\infty(\mathbb{R}^2)$ is the set of all functions that are in $\mathcal{C}^\infty(\IR^2)$, quasi-periodic with parameter $\alpha_0$, and vanish for large $|x_2|$; $\mathcal{C}_{\alpha_0}^\infty(\Omega)$ is the space of their restrictions to $\Omega$. 
The space $H^1_{\alpha_0}(\Omega)$ is the closure of $\mathcal{C}_{\alpha_0}^\infty(\Omega)$ in $H^1(\Omega)$.
This is a Hilbert space with the usual $H^1(\Omega)$ norm and inner product.

Every $u \in \mathcal{C}_{\alpha_0}^\infty(\Omega)$ may be represented as a Fourier series \cite[Prop.~2.6]{Pinto2020}:
\begin{equation} \label{eq:FourierExp}
u(x_1, x_2) = \sum_{n \in \mathbb{Z}} u_n(x_2) \re^{\ri \alpha_n x_1},
\qquad\text{where }\alpha_n := \alpha_0 + \frac{2\pi n}{L}, \; \text{for }n\in \mathbb{Z},
\end{equation}
$x_1\in[0,L]$, $x_2$ is such that $([0,L]\times\{x_2\})\cap D=\emptyset$,
and the coefficients $u_n$ are defined as
\begin{equation} \label{eq:FourierCoef}
u_n(x_2) := \frac{1}{L} \int_0^{L} \re^{-\ri\alpha_n x_1} u(x_1,x_2) \, \rd x_1 \qquad \text{for} \; n \in \mathbb{Z}.
\end{equation}

We introduce the quasi-periodic fractional Sobolev spaces on $\GHpm$:
\begin{equation*}
H^{\pm 1/2}_{\alpha_0}(\GHpm) := \bigg\{ v(x_1) = \sum_{n \in \mathbb{Z}} v_n \re^{\ri \alpha_n x_1} \; \biggl| \; \sum_{n \in \mathbb{Z}} (1+\alpha_n^2)^{\pm 1/2} \, |v_n|^2 < \infty \bigg\}.
\end{equation*}
It can be proved that $H^{1/2}_{\alpha_0}(\GHpm)$ is the space of the traces on $\GHpm$ of all functions of $H^1_{\alpha_0}(\Omega)$; see \cite[Lemma~2.27]{Pinto2020}.
The duality product between $H^{1/2}_{\alpha_0}(\GHpm)$ and $H^{-1/2}_{\alpha_0}(\GHpm)$ is given by 
\begin{equation*}
\langle  w , v \rangle_{\alpha_0, \GHpm} = L \sum_{n \in \mathbb{Z}}  w_n \overline{v_n},\qquad\text{for }  w(x_1)= \sum_{n \in \mathbb{Z}}  w_n \re^{\ri \alpha_n x_1}, \quad v(x_1)= \sum_{n \in \mathbb{Z}} v_n \re^{\ri \alpha_n x_1}.
\end{equation*}

We now specify a radiation condition on $\GHpm$, as we want that the scattered field propagates upward on $\GHp$ and the total field propagates downward on $\GHm$.

Recalling the coefficients $\alpha_n$ in \eqref{eq:FourierExp}, we first define the values
\begin{equation}\label{eq:BetaEps}
\beta_n[\varepsilon]:=\sqrt{k^2\varepsilon-\alpha_n^2}, \qquad n\in\IZ, \qquad 
\Re\varepsilon>0, \quad \Im\varepsilon\ge0.
\end{equation}
When $\varepsilon\notin\IR$ we select the complex square root with positive imaginary part, i.e.\ $\Im\beta_n[\varepsilon]>0$.
When $\varepsilon\in\IR$, then $\beta_n[\varepsilon]$ in \eqref{eq:BetaEps} is real and positive for all $n$ such that $\alpha_n^2\le k^2\varepsilon$, while $\beta_n[\varepsilon]=\ri\sqrt{\alpha_n^2-k^2\varepsilon}$ is purely imaginary with $\Im\beta_n[\varepsilon]>0$ for all other values of $n$.
In this case, $\beta_n[\varepsilon]$ is real for a finite number of indices $n$ (the propagative modes).
The remaining indices $n$ are the evanescent modes.
In particular, $\beta_0(\varepsilon^+)=-\kappa^+\sin\theta$ is real and positive.

\begin{rem}[Decomposition of quasi-periodic waves]\label{rem:Decomposition}
For all $n\in\IZ$ and $\varepsilon$ as in \eqref{eq:BetaEps}, we have $\alpha_n^2+(\beta_n[\varepsilon])^2=k^2\varepsilon$.
Let $w$ be any quasi-periodic solution of the Helmholtz equation $\Delta w+k^2\varepsilon w=0$ in a rectangle $R=(0,L)\times(x_2^\flat,x_2^\sharp)$ with constant $\varepsilon$.
Then, from separation of variables, there is a unique decomposition 
$$
w(\bx)=
\sum_{n\in\IZ} c_n^\uparrow\re^{\ri\alpha_n x_1}\re^{\ri\beta_n[\varepsilon] x_2}+
\sum_{n\in\IZ} c_n^\downarrow\re^{\ri\alpha_n x_1}\re^{-\ri\beta_n[\varepsilon] x_2}
=:w^\uparrow(\bx)+w^\downarrow(\bx), \qquad c_n^\downarrow,c_n^\uparrow\in\IC.$$
The fields $w^\uparrow$ and $w^\downarrow$ are the upward and downward propagating components of $w$.
The DtN operators allow to select one of these two components.
\end{rem}

We first focus on $\GHp$. 
In the semi-infinte strip $\Strip^+$ we write the total field $u$ as the sum of the incident field \eqref{eq:Incident} and a scattered field:
\begin{equation*}
u = u^\inc+u^\scat.
\end{equation*}
Since both $u$ and $u^\inc$ are smooth, quasi-periodic solutions of the Helmholtz equation in $\Strip^+$, also their difference $u^\scat$ enjoys the same properties and propagates upwards, thus
\begin{equation} \label{eq:Representation}
u^\scat(\bx)=\sum_{n\in\mathbb{Z}} u^\scat_n(H^+) 
\re^{\ri\alpha_n x_1}\re^{\ri\beta_n[\varepsilon^+] (x_2-H^+)},
\qquad \bx\in\Strip^+ \cup \GHp.
\end{equation}
To enforce the relation between the value \eqref{eq:Representation} of $u^\scat$ and that of its normal derivative in the formulation of the scattering problem, we define the Dirichlet-to-Neumann (DtN) operator $\DtN_{\varepsilon^+}^+$ as
\begin{align} \label{eq:DtN}
& \DtN_{\varepsilon^+}^+:H^{1/2}_{\alpha_0}(\GHp) \to H^{-1/2}_{\alpha_0}(\GHp), \\ \nonumber
& (\DtN_{\varepsilon^+}^+\phi)(x_1) := \ri \sum_{n\in\mathbb{Z}} \phi_n \beta_n[\varepsilon^+]\re^{\ri\alpha_nx_1}, 
\hspace{0.5cm} \text{for} \hspace{0.2cm} 
\phi(x_1)=\sum_{n\in\IZ}\phi_n \re^{\ri\alpha_nx_1} \in H^{1/2}_{\alpha_0}(\GHp).
\end{align}
Deriving with respect to $x_2$ the formula in Remark~\ref{rem:Decomposition} and comparing against the action \eqref{eq:DtN} of the DtN operator, we observe
that a Helmholtz solution $u^\scat$ in $\Strip^+$ propagates upwards
if $\gamma_{H^+}(\dxt u^\scat)=\DtN_{\varepsilon^+}^+(\gamma_{H^+} u^\scat)$, where $\gamma_{H^+}$ is the trace on $\GHp$. 
This is equivalent to the requirement that
$u^\scat$ admits the expansion \eqref{eq:Representation}.
If these equivalent conditions are met, we say that $u^\scat$ satisfies the radiation condition in $\Strip^+$.

Similarly, on $\GHm$, we require the total field $u$ to propagate downwards:
\begin{equation}\label{eq:ReprDown}
u(\bx)=\sum_{n\in\mathbb{Z}} u_n(H^-)\re^{-\ri\beta_n[\varepsilon^-] (x_2-H^-)} \re^{\ri\alpha_n x_1},    \qquad \bx\in\Strip^-\cup\GHm.
\end{equation}
The DtN operator $\DtN_{\varepsilon^-}^-$ on $\GHm$ is:
\begin{align} \label{eq:DtNmeno}
& \DtN_{\varepsilon^-}^-:H^{1/2}_{\alpha_0}(\GHm) \to H^{-1/2}_{\alpha_0}(\GHm), \\ \nonumber
& (\DtN_{\varepsilon^-}^- \phi)(x_1) := \ri \sum_{n\in\mathbb{Z}} \phi_n \beta_n[\varepsilon^-] \re^{\ri\alpha_n x_1}, \hspace{0.5cm} \text{for} \hspace{0.2cm} 
\phi(x_1)=\sum_{n\in\IZ}\phi_n  \re^{\ri\alpha_nx_1} \in H^{1/2}_{\alpha_0}(\GHm).
\end{align}
The downward-propagation condition on $u$ is equivalent to the requirement that $\gamma_{H^-}(-\dxt u)=\DtN_{\varepsilon^-}^-(\gamma_{H^-} u)$. The parameters $\alpha_n=k\sqrt{\varepsilon^+}\cos\theta+\frac{2\pi n}L$ enter the definition of the DtN operator $\DtN_{\varepsilon^-}^-$ on the lower boundary $\GHm$, but they depend on the value $\varepsilon^+$ of the material parameter $\varepsilon$ in the upper region $\Strip^+$, while $\varepsilon^-$ enters $\DtN_{\varepsilon^-}^-$ via the $\beta_n[\varepsilon^-]$ coefficients \eqref{eq:BetaEps}.

The operators $\DtN_{\varepsilon^\pm}^\pm:H^{1/2}_{\alpha_0}(\GHpm) \to H^{-1/2}_{\alpha_0}(\GHpm)$ are continuous \cite[Lemma 3.7]{Pinto2020}.
See \cite{Pinto2020} for more details on the DtN operators properties and their derivation.

Besides the permittivity $\varepsilon^\pm$, the operators $\DtN_{\varepsilon^\pm}^\pm$ \eqref{eq:DtN} and \eqref{eq:DtNmeno} also depend on $L,k,\alpha_0$, through the coefficients $\beta_n[\varepsilon^\pm]$ in \eqref{eq:BetaEps}; however, in the notation we only stress the dependence on $\varepsilon^\pm$ since the T-matrix method will involve DtN operators with differente permittivity values (see \S\ref{s:BVPlayer}).

In practical computations, we need to truncate the infinite series in the definitions of the DtN operators (\ref{eq:DtN}) and (\ref{eq:DtNmeno}), so we introduce the following operators.
\begin{defin}[Truncated DtN operator]\label{def:TruncatedDtN}
For $M\in\IN$, we set
\begin{equation} \label{eq:DtNtrunc}
(\DtN_{\varepsilon^\pm,M}^\pm \phi)(x_1) := 
\ri \sum_{n=-M}^M \phi_n \beta_n[\varepsilon^\pm] \re^{\ri\alpha_nx_1},
\qquad\forall \phi(x_1)=\sum_{n\in\IZ}\phi_n\re^{\ri\alpha_nx_1} \in H^{1/2}_{\alpha_0}(\GHpm).
\end{equation}
\end{defin}

The choice of the truncation order $M$ is crucial for the well-posedness of the BVP.
It can be shown \cite[Remark~2.6]{Monforte2025} that if 
\begin{align}\label{eq:MLarge}
\varepsilon^-\in\IR, \qquad 
M\ge{M_\star:=} \frac L{2\pi}\Big(\max\{\kappa^+,\kappa^-\} +|\alpha_0|\Big),
\end{align}
then 
$$
\Im \int_\GHpm\DtN_{\varepsilon^\pm,M}^\pm w\; \cw\rd s=\Im \int_\GHpm\DtN_{\varepsilon^\pm}^\pm w\; \cw\rd s
\qquad \forall w\in H^1_{\alpha_0,0}(\Omega).
$$
This identity is key to ensure the well-posedness of the TDG scheme \eqref{eq:TruncDG}, see \cite[Prop.~4.2]{Monforte2025}.

\subsection{Boundary value problem}\label{s:BVP}
Using the DtN operators $\DtN_{\varepsilon^\pm}^\pm$, we impose the radiation condition on $\GHpm$, obtaining the following BVP: given the incident wave 
$u^\inc(\bx) = \re^{\ri \kappa^+(x_1\cos \theta + x_2\sin \theta)}= \re^{\ri\alpha_0 x_1 - \ri \beta_0^+ x_2}$,
with $\theta \in (-\pi, 0)$, find $u \in H^1_{\alpha_0}(\Omega)$, with $\alpha_0 = \kappa^+ \cos \theta$, such that
\begin{equation} \label{eq:BVP}
\begin{cases}
\Delta u + k^2\varepsilon u = 0 & \hspace{0.5cm} \text{in} \; \Omega, \\
u = 0 & \hspace{0.5cm} \text{on} \; \GD, \\
\dn(u-u^\inc) - \DtN_{\varepsilon^+}^+(u-u^\inc) = 0 & \hspace{0.5cm} \text{on} \; \GHp, \\
\dn u - \DtN_{\varepsilon^-}^-u = 0 & \hspace{0.5cm} \text{on} \; \GHm.
\end{cases}
\end{equation}

Numerically, we approximate a truncated BVP, making use of the truncated DtN operators $\DtN_{\varepsilon^\pm,M}^\pm$ \eqref{eq:DtNtrunc}, finding a solution $u^M \in H^1_{\alpha_0}(\Omega)$, for $M \in \IN$.

\begin{rem}[Well-posedness of the BVP] \label{rem:well-posed}
The BVP \eqref{eq:BVP} is well-posed for every wavenumber $k$ under non-trapping assumptions \cite[Theorem~3.5]{BonnetBD1994}, \cite[Theorem~3.3]{Monforte2025}: geometrically, they mean that a point moving along a vertical half line from any $(x_1,0)$, either upwards or downwards, does not enter any impenetrable obstacle, and that at every material interface the value of $\varepsilon$ increases. 
If these conditions are not satisfied, \eqref{eq:BVP} is well-posed in its weak form for every value of $k$ except possibly for an increasing sequence of ``singular frequencies'' $(k_m)_{m \geq 1}$ that tends to infinity with~$m$ \cite[Theorem 3.2--3.4]{BonnetBD1994}. The anomalous cases $k=k_m$ correspond to situations where the system allows a solution known as ``guided mode'' \cite[\S4.1]{BonnetBD1994}, ``Rayleigh--Bloch wave in the radiation continuum'', or ``bound state in the continuum'' \cite{Marinica2008,hsu_bound_2016}.
The truncated problem with $M\ge M_\star$ as in \eqref{eq:MLarge} is well-posed under the same assumptions, see \cite[Corollary~3.4]{Monforte2025}.
\end{rem}

\section{Scattering problems on layered structures}\label{s:LayeredStructures}
\subsection{Periodic layered structures}\label{s:Layers}
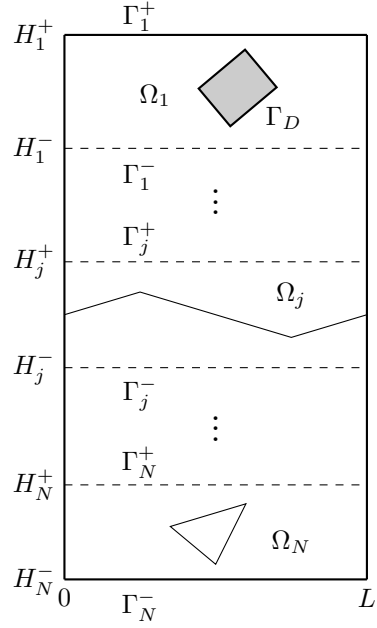
\begin{figure}[tbp]
\centering
\begin{tikzpicture}
\draw[thick] (0,-2.5) -- (0,4.7);
\draw[thick] (4,-2.5) -- (4,4.7);
\draw[thick] (0,4.7) -- (4,4.7);
\draw[thick] (0,-2.5) -- (4,-2.5);
\draw[dashed] (0,3.2) -- (4,3.2);
\draw[dashed] (0,1.7) -- (4,1.7);
\draw[dashed] (0,0.3) -- (4,0.3);
\draw[dashed] (0,-1.25) -- (4,-1.25);
\fill[gray!40!white, draw=black, thick, rotate around={40:(2.7,3.2)}] (2.5,3.75) rectangle (3.3,4.4);
\draw (0,1) -- (1,1.3) -- (3,0.7) -- (4,1);
\coordinate (A) at (1.4,-1.8);
\coordinate (B) at (2,-2.3);
\coordinate (C) at (2.4,-1.5);
\draw (A) -- (B) -- (C) -- cycle;
\draw (0,4.7) node[left]{$H_1^+$};
\draw (0,3.2) node[left]{$H_1^-$};
\draw (0,1.7) node[left]{$H_j^+$};
\draw (0,0.3) node[left]{$H_j^-$};
\draw (0,-1.25) node[left]{$H_N^+$};
\draw (0,-2.5) node[left]{$H_N^-$};
\node[] at (1.2,3.9) {$\Omega_1$};
\node[] at (3,1.25) {$\Omega_j$};
\node[] at (3,-2) {$\Omega_N$};
\draw (0,-2.5) node[below]{$0$};
\draw (4,-2.5) node[below]{$L$};
\node[] at (1,4.95) {$\Gamma_1^+$};
\node[] at (1,2.85) {$\Gamma_1^-$};
\node[] at (1,2) {$\Gjp$};
\node[] at (1,-0.05) {$\Gjm$};
\node[] at (1,-0.95) {$\Gamma_N^+$};
\node[] at (1,-2.85) {$\Gamma_N^-$};
\node[] at (2.9,3.6) {$\GD$};
\node[] at (2,2.6) {\Large$\vdots$};
\node[] at (2,-0.4) {\Large$\vdots$};
\end{tikzpicture}  
\caption{The layered truncated domain $\Omega$. Continuous lines represent the material interfaces, dashed horizontal lines separate the domain layers.}
\label{fig:layers}
\end{figure}

We now focus on periodic domains $\Omega_0$ (and consequently their truncations $\Omega$) composed of $N$ stacked non-overlapping regions, possibly at some distance from one another.
For $N\in \IN$ and the ``heights'' $$H_N^-<H_N^+\le H_{N-1}^-<H_{N-1}^+\le\cdots\le H_2^- <H_2^+\le H_1^-<H_1^+,$$
we define $N$ layers $\Omega_j$, their horizontal boundaries $\Gjpm$, and the adjacent unbounded upper/lower regions by
\begin{align} 
\nonumber
\Omega_j := &\{\bx\in\Omega_0:\; 0<x_1<L,\; H_j^-<x_2<H_j^+\}, \\
\Gjpm:=&\{\bx\in\IR^2: 0\le x_1\le L,\; x_2=H_j^\pm \}, \nonumber
\\
\Strip_j^\pm:=&\{\bx\in\IR^2: 0 < x_1 < L,\; x_2\gtrless H_j^\pm \},
\qquad\qquad j = 1, \ldots, N,
\label{eq:Omega_j}
\end{align}
and define $\Omega$ as the union of all the layers and the regions between them:
\begin{equation}\label{eq:OmegaLayers}
\Omega := \bigcup_{j=1}^N \Omega_j \cup \bigcup_{j=1}^{N-1} (0,L) \times [H_{j+1}^+,H_j^-]
=\big((0,L)\times(H_N^-,H_1^+)\big)\cap\Omega_0.
\end{equation} 
For $j=1,\dots,N$, we denote by $\varepsilon_j^{+}$ and $\varepsilon_j^{-}$ the values of the relative permittivity $\varepsilon$ immediately above $\Gamma_j^+$ and immediately below $\Gamma_j^-$, respectively.
We assume that the PEC obstacle $D\cap\{0<x_1<L\}$ and all discontinuities in the permittivity $\varepsilon$ are contained in $\bigcup_{j=1}^N\Omega_j$.
This means that $\varepsilon(\bx)=\varepsilon_1^+$ for $x_2\ge H_1^+$, $\varepsilon(\bx)=\varepsilon_N^-$ for $x_2\le H_N^-$, and that between two subsequent regions $\Omega_j$ and $\Omega_{j+1}$ the relative permittivity $\varepsilon$ assumes a constant value:
\begin{align*}
\varepsilon(\bx) & = \varepsilon_j^- = \varepsilon_{j+1}^+ \qquad \text{in } \{ H_{j+1}^+ \leq x_2 \leq H_j^- \} ,
\qquad 1\le j\le N-1.
\end{align*}
Setting the quantities introduced in \S\ref{s:Domain} as
$$
H^+ = H_1^+, \quad H^- = H_N^-, \quad
\GHp=\Gamma_1^+, \quad \GHm=\Gamma_N^-, \quad
\varepsilon^+ = \varepsilon_1^+, \quad\varepsilon^- = \varepsilon_N^-, 
$$
and with the incoming plane wave $u^\inc$ as in \S\ref{s:BVP}, we consider the BVP \eqref{eq:BVP} in the special case where the domain $\Omega$ is layered as in \eqref{eq:OmegaLayers}.

An example of layered domain is depicted in Figure \ref{fig:layers}: the continuous lines, periodic in $x_1$, represent the interfaces separating different materials, i.e.\ regions with constant $\varepsilon$. Dashed horizontal lines separate the domain layers $\Omega_j$.

In the following, we choose not to directly solve the scattering problem \eqref{eq:BVP} on the global domain $\Omega$, but to solve instead several smaller problems on every layer $\Omega_j$ and then couple them to obtain the solution of \eqref{eq:BVP}.
To do this efficiently, we adapt the T-matrix method \cite{Waterman1965}, which has been applied to multiple scattering problems \cite{Hawkins2024,Martin2006,Monforte2026,Mishchenko1996}.
We derive a T-matrix for each layer and use a translation operator to couple subsequent layers.

\subsection{Local quasi-periodic BVP on a single layer}\label{s:BVPlayer}
The T-matrix method requires solving several scattering problems for each layer.
We consider BVPs posed on $\Omega_j$ with homogeneous media in the complement $\Strip_j^+\cup\Strip_j^-$ setting 
$\varepsilon=\varepsilon_j^\pm$ on $\Strip_j^\pm$.
This means that we extend the values of $\varepsilon$ from the upper and lower faces $\Gjpm$ of $\Omega_j$ outwards as a constant.

\begin{figure}[htbp]
\centering
\begin{tikzpicture}
\draw (0,3) -- (4,3);
\draw (0,1.5) -- (4,1.5);
\draw (0,0) -- (4,0);
\draw (0,0) -- (0,3);
\draw (4,0) -- (4,3);
\node[] at (2,0.5) {$\varepsilon_j^-$};
\node[] at (2,2.5) {$\varepsilon_j^+$};
\node[] at (1,3.25) {$\Gjp$};
\node[] at (1,-0.3) {$\Gjm$};
\fill[gray!40!white, draw=black, thick, rotate around={45:(2.5,0.95)}] (2.35,0.9) rectangle (3.35,1.9);
\draw[-{Stealth[slant=0]}, thick, rotate around={45:(2.5,0.95)}](3.35,1.4)--(2.85,1.4);
\node[] at (2.3,1.7) {$\mathbf{n}$};
\draw (0,0) node[left]{$H_j^-$};
\draw (0,3) node[left]{$H_j^+$};
\draw (0,0) node[below]{$0$};
\draw (4,0) node[below]{$L$};
\draw [dashed,->](0,1.5) -- (4.7,1.5);
\draw [dashed,->](0,3) -- (0,3.7);
\node[] at (4.7,1.2) {$x_1$};
\node[] at (-0.3,3.7) {$x_2$};
\draw [-{Stealth[slant=0]}, thick](2,3) -- (2,3.6);
\node[] at (2.2,3.3) {$\mathbf{n}$};
\draw [-{Stealth[slant=0]}, thick](2,0) -- (2,-0.6);
\node[] at (2.2,-0.3) {$\mathbf{n}$};
\node[] at (3.2,-0.5) {$\Strip_j^-$};
\node[] at (3.2,3.4) {$\Strip_j^+$};
\end{tikzpicture}
\caption{A sketch of the layer $\Omega_j$.}
\end{figure}
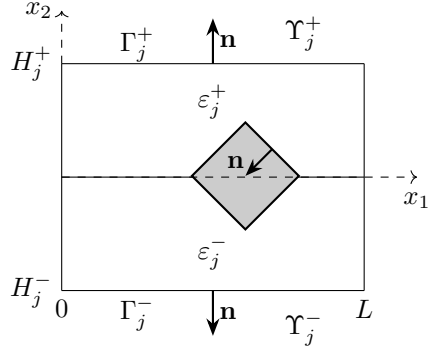

The incident wave $u_{\inc,j}$ is a quasi-periodic function with parameter $\alpha_0 = \kappa_1^+ \cos \theta$, inherited from the global scattering problem on $\Omega$, even if $\varepsilon_j^+\neq \varepsilon_1^+$ for $j>1$.
Since the incident wave on $\Omega_j$ is the result of the wave scattered by other layers above or below $\Omega_j$, and all of these are quasi-periodic functions, in general $u_{\inc,j}$ is the combination of an upward-propagating wave $u_{\inc,j}^\uparrow$ coming from below and a downward-propagating wave $u_{\inc,j}^\downarrow$ coming from above.

Following \eqref{eq:BetaEps}, for every layer $\Omega_j$, the coefficients $\beta_{j,n}^\pm$ are defined as 
\begin{equation}\label{eq:BetaJN}
\beta_{j,n}^\pm:=\beta_n[\varepsilon_j^\pm]=\sqrt{k^2\varepsilon_j^\pm-\alpha_n^2},
\qquad j\in\{1,\ldots,N\}, \qquad n\in \IZ.
\end{equation}
Differently from \S\ref{s:periodic_general} and the global BVP \eqref{eq:BVP}, here we allow $\varepsilon_j^+\notin\IR$, i.e.\ the upper region $\Strip_j^+$ may contain an absorbing medium.

We impose a radiation condition requiring that the horizontal boundaries $\Gjpm$ are transparent to waves propagating away from the grating.
This means that in $\Gjp \cup \Strip_j^+$ the difference $u_j-u_{\inc,j}^\downarrow$ propagates upward and in $\Gjm \cup \Strip_j^-$ the difference $u_j-u_{\inc,j}^\uparrow$ propagates downward, the BVP solution $u_j$ being the total field.
Using the DtN operators $\DtN_{\varepsilon_j^\pm}^\pm$ \eqref{eq:DtN}--\eqref{eq:DtNmeno}, we obtain the following BVP: given the incident waves
$u_{\inc,j}^\downarrow\in H^{1/2}_{\alpha_0}(\Gjp)$, $u_{\inc,j}^\uparrow\in H^{1/2}_{\alpha_0}(\Gjm)$, 
find $u_j\in H^1_{\alpha_0}(\Omega_j)$ such that
\begin{equation} \label{eq:bvp_single_layer}
\begin{cases}
\Delta u_j + k^2\varepsilon u_j = 0 & \hspace{0.5cm} \text{in} \; \Omega_j, \\
u_j = 0 & \hspace{0.5cm} \text{on} \; \Gamma_D\cap\partial\Omega_j, \\
\dn(u_j-u_{\inc,j}^\downarrow) - \DtN_{\varepsilon_j^+}^+(u_j-u_{\inc,j}^\downarrow) = 0 & \hspace{0.5cm} \text{on} \; \Gjp, \\
\dn(u_j-u_{\inc,j}^\uparrow) - \DtN_{\varepsilon_j^-}^-(u_j-u_{\inc,j}^\uparrow) = 0 & \hspace{0.5cm} \text{on} \; \Gjm.
\end{cases}
\end{equation}
We can also define a truncated BVP replacing the DtN operators $\DtN_{\varepsilon_j^\pm}^\pm$ with their truncation $\DtN_{\varepsilon_j^\pm,M}^\pm$.
These are the BVPs that we approximate numerically.

\begin{rem}[Transmitted and reflected waves]
In $\Gjp \cup \Strip_j^+$, using the Fourier modal expansion \eqref{eq:FourierExp} of $u_j$ on $\Gjp$, we define the upward-propagating wave $u_j^\uparrow$ as
\begin{equation} \label{eq:representation_up_layer}
\begin{split}
u_j^\uparrow(\bx) := u_j(\bx)-u_{\inc,j}^\downarrow(\bx)
& = u_{j,T}^\uparrow(\bx) + u_{j,R}^\uparrow(\bx) \\
& = \sum_{n\in\mathbb{Z}} \left[u_{T,n}^\uparrow(H_j^+) + u_{R,n}^\uparrow(H_j^+)\right] \re^{\ri\beta_{j,n}^+ (x_2-H_j^+)} \re^{\ri\alpha_n x_1}, 
\end{split}
\end{equation}
where the subscripts $_T$ and $_R$ tag the ``transmitted'' and the ``reflected'' components, i.e.\ the solutions of BVP \eqref{eq:bvp_single_layer} with $u_{\inc,j}^\downarrow=0$ and $u_{\inc,j}^\uparrow=0$, respectively.
Similarly, on $\Gjm \cup \Strip_j^-$, the downward-propagating wave $u_j^\downarrow$ is defined as
\begin{equation} \label{eq:representation_down_layer}
\begin{split}
u_j^\downarrow(\bx) := u_j(\bx)-u^\uparrow_{\inc,j}(\bx)
& = u_{j,T}^\downarrow(\bx) + u_{j,R}^\downarrow(\bx) \\
& =\sum_{n\in\mathbb{Z}} \left[u_{T,n}^\downarrow(H_j^-) + u_{R,n}^\downarrow(H_j^-)\right] \re^{-\ri\beta_{j,n}^- (x_2-H_j^-)} \re^{\ri\alpha_n x_1},
\end{split}
\end{equation}
where in this region the transmitted and the reflected components $u_T$ and $u_R$ solve BVP \eqref{eq:bvp_single_layer} with $u_{\inc,j}^\uparrow=0$ and $u_{\inc,j}^\downarrow=0$, respectively.
\end{rem}

\section{The T-matrix method}\label{s:Tmat}

The T-matrix associated to a given scatterer is an infinite matrix that maps the incoming field to the field scattered by the obstacle.
Both the incoming and the scattered field are represented by their coefficient expansion in appropriate Fourier bases.

We first define the T-matrix associated with the BVP \eqref{eq:bvp_single_layer} posed on a single periodic layer $\Omega_j$ in \S\ref{s:TmatSingle}.
Then in \S\ref{s:TmatCoupling} we show how to couple the T-matrices associated to several layers $\Omega_j$ to solve the BVP \eqref{eq:BVP} posed on the layered structure \eqref{eq:OmegaLayers}.
For practical computations, the infinite T-matrix has to be truncated to a finite one (\S\ref{s:TmatApprox}) and its entries approximated by a numerical method: in \S\ref{s:TDG} we describe how to assemble it using the DtN-TDG method.

\subsection{The T-matrix for a single layer}\label{s:TmatSingle}
We use the expansions \eqref{eq:representation_up_layer} and \eqref{eq:representation_down_layer} to derive a relation between the incident and the reflected/transmitted field for BVP \eqref{eq:bvp_single_layer}.
We define the following basis functions, for $n \in \IZ$:
\begin{equation} \label{eq:psi_phi}
\begin{aligned}
\Psidown{j,n}(\bx) & := \re^{\ri\alpha_n x_1} \re^{-\ri\beta_{j,n}^+(x_2 - H_j^+)}, \qquad
&&\bx\in \Strip_j^+ \cup \Gjp, \\
\Phiup{j,n}(\bx) & := \re^{\ri\alpha_n x_1} \re^{\ri\beta_{j,n}^+(x_2 - H_j^+)}, &&\bx\in \Strip_j^+ \cup \Gjp,\\
\Psiup{j,n}(\bx) & := \re^{\ri\alpha_n x_1} \re^{\ri\beta_{j,n}^-(x_2 - H_j^-)}, &&\bx\in \Strip_j^- \cup \Gjm,\\
\Phidown{j,n}(\bx) & := \re^{\ri\alpha_n x_1} \re^{-\ri\beta_{j,n}^-(x_2 - H_j^-)}, &&\bx\in \Strip_j^- \cup \Gjm.
\end{aligned}
\end{equation}
The letters $\Psi$ and $\Phi$ denote the incoming and the outgoing basis functions, respectively.
The arrows ${}^{\uparrow\downarrow}$ denote the directions of propagation.
We use the basis functions \eqref{eq:psi_phi} to expand the incident and reflected/transmitted field $u_j$, solution of \eqref{eq:bvp_single_layer}, in $\Strip_j^\pm$:
\begin{equation} \label{eq:basis_expansion}
\begin{aligned} 
u_{\inc,j}^\downarrow(\bx) & = \sum_{n\in\IZ} a_{j,n}^\downarrow \Psidown{j,n}(\bx), \qquad 
& u_{j,T}^\uparrow(\bx) + u_{j,R}^\uparrow(\bx) & = \sum_{n\in\IZ} b_{j,n}^\uparrow \Phiup{j,n}(\bx), \qquad \bx\in \Strip_j^+ \cup \Gjp,\\
u_{\inc,j}^\uparrow(\bx) & = \sum_{n\in\IZ} a_{j,n}^\uparrow \Psiup{j,n}(\bx), 
& u_{j,R}^\downarrow(\bx) + u_{j,T}^\downarrow(\bx) & = \sum_{n\in\IZ} b_{j,n}^\downarrow \Phidown{j,n}(\bx), \qquad \bx\in \Strip_j^- \cup \Gjm.
\end{aligned}
\end{equation}
We denote by $\adown_j=(a_{j,n}^\downarrow)_{n\in\IZ}$ the Fourier expansion coefficient (infinite) vector of $u_{\inc,j}^\downarrow$ and by $\aup_j=(a_{j,n}^\uparrow)_{n\in\IZ}$ the vector associated to $u_{\inc,j}^\uparrow$.
Similarly, $\bup=(b_{j,n}^\uparrow)_{n\in\IZ}$ and $\bdown=(b_{j,n}^\downarrow)_{n\in\IZ}$ are the coefficient vectors of the outgoing field propagating upward in $\Strip_j^+$ and downward in $\Strip_j^-$, respectively.
Since the Helmholtz equation is linear, the relation between the coefficients of the incident field and the coefficient of the reflected/transmitted field is linear and can be expressed through a matrix-vector product:
\begin{align*}
& \bup_j = \Tup_j \aup_j + \Rup_j \adown_j, \\
& \bdown_j = \Rdown_j \aup_j + \Tdown_j \adown_j,
\end{align*}
where $\bT_j^\uparrow$, $\bT_j^{\downarrow}$ are the infinite matrices that map the incident upward/downward propagating wave into the transmitted field, and $\mathbf{R}_j^\uparrow$, $\mathbf{R}_j^\downarrow$ are the matrices that map into the reflected field.
We relate the incident vector to the reflected/transmitted vector in the following way:
\begin{equation} \label{eq:tmat}
\begin{pmatrix} \bup_j \\ \bdown_j \end{pmatrix} = \Tmat_j \begin{pmatrix} \aup_j \\ \adown_j \end{pmatrix}, \qquad \Tmat_j := \begin{pmatrix} \Tup_j & \Rup_j \\ \Rdown_j & \Tdown_j \end{pmatrix}.
\end{equation}
We call the infinite $2\times2$-block matrix $\Tmat_j$ the T-matrix associated to the layer $\Omega_j$.
The entries of $\Tmat_j$ depend only on the layer $\Omega_j$, i.e.\ its geometry, the height $H_j^+-H_j^-$, and the values of $\varepsilon$ inside $\Omega_j$, but they do not depend on the particular incident field.

This means that, once we have computed (or approximated) $\Tmat_j$, for any quasi-periodic incoming field $u_{\inc,j}$ in \eqref{eq:bvp_single_layer} we can compute the corresponding scattered field on $\Strip_j^\pm$ by multiplying $\Tmat_j$ with the coefficient vector $(\aup_j,\adown_j)$ of $u_{\inc,j}$ and using the expansion \eqref{eq:basis_expansion}. 
Most importantly, we are able to couple layers easily, as we show in \S\ref{s:TmatCoupling}.

\begin{rem}[T-matrix symmetry: translation]\label{rem:Translate}
If two layers $\Omega_j$, $\Omega_{j_*}$ are vertical translations of one another (i.e.\ $\bx\in\Omega_j \iff \bx+(0,z)\in \Omega_{j_*}$ and $\varepsilon(\bx)=\varepsilon_*(\bx+(0,z))$ for $\bx\in\Omega_j$ and $z\in\IR$), then the corresponding T-matrices coincide: $\Tmat_j=\Tmat_{j_*}$.
The corresponding basis functions in \eqref{eq:psi_phi} are translations of one another, e.g.\ $\Psidown{j,n}(\bx)=\Psidown{j_*,n}(\bx+(0,z))$.
\end{rem}

\begin{rem}[T-matrix symmetry: vertical reflection]\label{rem:Reflect}
	Let $\mho_j$ be the reflection of $\Omega_j$ with respect to the horizontal axis $x_2 = 0$, i.e.\
	$(x_1, x_2)\in\mho_j \iff (x_1, -x_2)\in\Omega_j$,
	with relative permittivity $\widetilde\varepsilon(x_1,x_2)=\varepsilon(x_1,-x_2)$.
	Then the T-matrix $\widetilde{\Tmat}_j$ associated to $\mho_j$ is given by
	\begin{equation*}
		\widetilde{\Tmat}_j = \begin{pmatrix} \Tdown_j & \Rdown_j \\ \Rup_j & \Tup_j \end{pmatrix},
	\end{equation*}
	i.e.\ the diagonal blocks $\Tup_j \leftrightarrow \Tdown_j$ and the off-diagonal blocks 
	$\Rup_j \leftrightarrow \Rdown_j$ are swapped with respect to $\Tmat_j$ 
	(cf.~\eqref{eq:tmat}). This holds for any quasi-periodicity parameter $\alpha_0$, 
	since the reflection acts only on $x_2$ and preserves quasi-periodicity in $x_1$.
\end{rem}

\subsection{T-matrix method for a layered domain} \label{s:TmatCoupling}

We now derive the coupled system that allows to find the solution $u$ of the global BVP \eqref{eq:BVP} on the layered domain $\Omega$ in \eqref{eq:OmegaLayers}.
The idea of describing each layer by a ``layer matrix'' mapping incident to reflected and transmitted wave amplitudes, 
and of coupling the layers through ``propagation matrices'', is classical in the study of stratified media and arrays of slabs, 
including multilayer gratings; see e.g. \cite{MontielSquireBennetts2015} and the references therein. 
In the approach we propose, the layer matrices are computed using the DtN-TDG method described in the following section. 

To obtain input--output relations linking the wave amplitudes on the first layer with those on the last layer, two popular approaches exist.
In the first, a recursive relation relates the layer matrices for adjacent layers to the layer matrix for a pair of layers as a whole; this approach is often known as recursive S-matrix method or Redheffer star product \cite{Beutel21, Whittaker1999, Yoon21}. In the second, the layer matrices are assembled in a single block-banded matrix, as in \cite{Liu2026}. We adopt the latter approach, as detailed in the following.

We now consider $N$ stacked layers $\Omega_1, \ldots, \Omega_N$ with their respective T-matrices $\Tmat_1, \ldots, \Tmat_N$.
The datum $u^\inc$ is a quasi-periodic downward propagating incident wave on the first layer.
We denote by $\adown_1$ its expansion coefficient vector, which is assumed to be known. 
We first define the distance between subsequent domain $\Omega_{j-1}$ and $\Omega_j$ by $\Delta H_{j-1,j} := H_{j-1}^- - H_j^+ \ge 0$.
Using the expansion \eqref{eq:representation_down_layer}, we represent the downward propagating wave $u_{j-1}^\downarrow = u_{j-1,R}^\downarrow+u_{j-1,T}^\downarrow$ outgoing from $\Omega_{j-1}$ as an incident wave into the lower layer $\Omega_j$:
\begin{align*}
u_{j-1}^\downarrow(\bx) 
=&\sum_{n\in\IZ}b^\downarrow_{j-1,n}\Phidown{j-1,n}(\bx)\\
= & \sum_{n\in\IZ} b^\downarrow_{j-1,n}\re^{-\ri \beta_{j-1,n}^-(x_2-H_{j-1}^-)} \re^{\ri \alpha_n x_1} \\
= & \sum_{n\in\IZ} b^\downarrow_{j-1,n}\re^{\ri \beta_{j-1,n}^-(H_{j-1}^- - H_j^+)}\re^{-\ri \beta_{j-1,n}^-(x_2-H_j^+)} \re^{\ri \alpha_n x_1} \\
= & \sum_{n\in\IZ} b^\downarrow_{j-1,n}\re^{\ri \beta_{j,n}^+ \Delta H_{j-1,j} } \Psidown{j,n}(\bx),
\qquad \qquad \text{for }\bx\in (0, L) \times [ H_j^+, H_{j-1}^- ],
\end{align*}
where we used that $\beta_{j,n}^+ = \beta_{j-1,n}^-$ from \eqref{eq:BetaJN} and $\varepsilon_j^+=\varepsilon_{j-1}^-$.
Since this sum must be equal to the expansion of the downward incident wave $u_{\inc,j}^\downarrow$ into $\Gjp$ (recall \eqref{eq:basis_expansion}), 
$a^\downarrow_{j,n} = \re^{\ri \beta_{j,n}^+ \Delta H_{j-1,j}}\,b^\downarrow_{j-1,n}$ for $n\in\IZ$, and so
\begin{equation} \label{eq:coupling_down}
\adown_j = \bD_{j-1,j} \bdown_{j-1}, \qquad \bD_{j-1,j} := \mathrm{diag}(\re^{\ri \beta_{j,n}^+ \Delta H_{j-1,j}})_{n\in\IZ}.
\end{equation}
The matrix $\bD_{j-1,j}$ is a propagation matrix that takes into account the relative position of the layers.
If $H_{j-1}^-=H_j^+$, the two layers are adjacent and $\bD$ is the identity matrix.

Analogously, the upward propagating wave on $\Gjm$ can be written as the translation of the upward propagating scattered wave $u_{j+1}^\uparrow(\bx)$ coming from $\Omega_{j+1}$, leading to
\begin{equation} \label{eq:coupling_up}
\aup_j = \bD_{j,j+1} \bup_{j+1},
\end{equation}
where $\bD_{j,j+1}$ is the same diagonal matrix as in \eqref{eq:coupling_down} but associated to different layers.

Using the T-matrix relation \eqref{eq:tmat} and the relations \eqref{eq:coupling_down}--\eqref{eq:coupling_up} we get, for $j=2, \ldots, N-1$,
\begin{align*}
& \bup_j = \Tup_j \aup_j + \Rup_j \adown_j = \Tup_j \bD_{j,j+1} \bup_{j+1} + \Rup_j \bD_{j-1,j} \bdown_{j-1}, \\
& \bdown_j = \Rdown_j \aup_j + \Tdown_j \adown_j = \Rdown_j \bD_{j,j+1} \bup_{j+1} + \Tdown_j \bD_{j-1,j} \bdown_{j-1}.
\end{align*}
The analogous relations for the first and last layer are slightly different, since $\adown_1$ is known and cannot be written in terms of waves outgoing from other layers, and $\aup_N=0$ from the homogeneous boundary condition on $\GHm$ in \eqref{eq:BVP}:
\begin{align*}
& \bup_1 = \Tup_1 \aup_1 + \Rup_1 \adown_1 = \Tup_1 \bD_{1,2} \bup_2 + \Rup_1 \adown_1, \\
& \bdown_1 = \Rdown_1 \aup_1 + \Tdown_1 \adown_1 = \Rdown_1 \bD_{1,2} \bup_2 + \Tdown_1 \adown_1, \\
& \bup_N = \Tup_N \aup_N + \Rup_N \adown_N = \Rup_N \bD_{N-1,N} \bdown_{N-1}, \\
& \bdown_N = \Rdown_N \aup_N + \Tdown_N \adown_N = \Tdown_N \bD_{N-1,N} \bdown_{N-1}.
\end{align*}

We write the resulting linear system for $\bup_1, \bdown_1, \ldots, \bup_N, \bdown_N$:
\begin{equation} \label{eq:coupling_system}
\mathbb S \begin{bmatrix}
\bup_1 \\ \bdown_1 \\ \vdots \\ \bup_j \\ \bdown_j \\ \vdots \\ \bup_N \\ \bdown_N
\end{bmatrix}
=
\begin{bmatrix}
\Rup_1 \adown_1 \\
\Tdown_1 \adown_1 \\
0 \\ \\ \\
\vdots \\ \\ \\
0
\end{bmatrix},
\end{equation}
where
\begin{small}
\begin{equation*} \label{eq:coupling_matrix}
\mathbb S =\begin{bmatrix}
\,\BoxMat{\bI}{\bzero}{\bzero}{\bI} \,\, \BoxMat{-\Tup_1 \bD_{1,2}}{\bzero}{-\Rdown_1 \bD_{1,2}}{\bzero} & & \\
\qquad \qquad \qquad \qquad \ddots & \\
& \hspace*{-7pt}\BoxMat{\bzero}{-\Rup_j \bD_{j-1,j}}{\bzero}{-\Tdown_j \bD_{j-1,j}} \,\, \BoxMat{\bI}{\bzero}{\bzero}{\bI} \,\, \BoxMat{-\Tup_j \bD_{j,j+1}}{\bzero}{-\Rdown_j \bD_{j,j+1}}{\bzero} & \\
& \qquad \qquad \qquad \qquad \qquad \qquad \qquad \qquad \ddots &  \\
& &  \hspace*{-7pt}\BoxMat{\bzero}{-\Rup_N \bD_{N-1,N}}{\bzero}{-\Tdown_N \bD_{N-1,N}} \,\, \BoxMat{\bI}{\bzero}{\bzero}{\bI}\,
\end{bmatrix}.
\end{equation*}
\end{small}

Here $\bI$ and $\bzero$ denote the identity and null matrix respectively.
The matrix $\mathbb S$ is a banded block matrix.
To solve the coupled system \eqref{eq:coupling_system} we can use a Band Gaussian Elimination algorithm \cite[\S4.3.1]{Golub2013}, whose complexity grows only linearly in $N$.

A tridiagonal operator matrix related to $\mathbb S$ is used in \cite[eq.~(4.18)]{PintoAylwinJH2021} and \cite[eq.~(36)--(37)]{ChoBarnett2015} for a 2D quasi-periodic layered Helmholtz problem, where each layer comprises a single (piecewise) smooth interface separating two homogeneous regions, and each matrix block corresponds to a singular boundary integral operator or potential.

\begin{rem}[Repeated layers and re-usability]
If two or more layers are equal, by Remark~\ref{rem:Translate} their T-matrices coincide.
So each of the four blocks $\Tup,\Tdown,\Rup,\Rdown$ has to be computed only once reducing the computational cost of the linear system assembly.
\end{rem}

\begin{rem}[Incident wave coming from below]
The scheme \eqref{eq:coupling_system} refers to a downward-propagating incident wave coming from above, as stated in \eqref{eq:Incident}, but it could be employed in the case of an upward-propagating wave incident on $\Gamma_N^-$, as long as it is quasi-periodic with the same parameter $\alpha_0 = \kappa_1^+ \cos \theta$.
In this case, at the right hand side of \eqref{eq:coupling_system} the terms $\Tup_N \aup_N$ and $\Rdown_N \aup_N$ appear as the last two entries, where $\aup_N$ is the expansion vector of the upward-propagating incident field.
\end{rem}

\subsection{Numerical approximation of the T-matrix}\label{s:TmatApprox}
Each matrix $\Tmat_j$ in \eqref{eq:tmat}, and consequently each block of the matrix $\mathbb S$ in \eqref{eq:coupling_system}, is by definition infinite-dimensional.
In computations, $\Tmat_j$ can be cropped to $\Tmat_j^M$ by truncating the Fourier expansion of all waves on the interfaces $\Gjpm$ to $n\in\{-M,\ldots, M\}$, for some $M\in\IN$.

When each local BVP \eqref{eq:bvp_single_layer} is solved by truncating the DtN operators with $\DtN_{\varepsilon_j^\pm,M}^\pm$ defined in \eqref{eq:DtNtrunc}, we choose the same truncation order $M$ also for $\Tmat_j^{ M}$.
Thus, $\mathbf{a}_j^{\uparrow\downarrow}, \mathbf{b}_j^{\uparrow\downarrow} \in \IC^{2M+1}$ and $\Tmat_j^{ M} \in \IC^{(4M+2)\times(4M+2)}$ for all $j=1,\ldots,N$.
In practice, the truncation order $M$ is chosen in dependence on the wavenumber $k$ and the period $L$ of the scatterer, as discussed in Re\-mark~\ref{rem:truncation}.

To assemble $\Tmat_j^{ M}$, we select as incident fields the $4M+2$ basis functions $\Psidown{j,n}$ and $\Psiup{j,n}$ \eqref{eq:psi_phi}, for $n \in \{-M, \ldots, M\}$, then we solve numerically the corresponding scattering problems in the layer $\Omega_j$, and finally we derive the Fourier coefficients of the scattered wave.
For each incident field, these coefficients constitute a column of the matrix $\Tmat_j^{ M}$.

The solution vector of the (truncated) system \eqref{eq:coupling_system} determines the solution $u$ on the complement of the layers $\Omega_j$, i.e.\ on the semi-infinite strips $\Strip^\pm$ and on the gaps $(0,L)\times[H_{j+1}^+,H_j^-]$.
To compute the solution inside the layers $\Omega_j$ a further post-processing step is needed; see Remark~\ref{rem:PostProc} for a concrete recipe.

The convergence of the T-matrix scheme and its stability with respect to perturbations of the data are classical for obstacle scattering \cite{KristenssonRammStrom1983,GaneshHawkinsHiptmair2012}. 
A well-known source of instability in classical T-matrix algorithms is due to the use of null-field equations to compute the T-matrix itself, see \cite[p.~543]{Mishchenko1996}: since in \S\ref{s:TDG} we compute the matrices using the TDG scheme, this issue does not arise (see however the comment on the instabilities due to the propagative plane wave basis and possible remedies in \S\ref{s:TDG}).
Coupling the layers as in \eqref{eq:coupling_system} through outgoing amplitudes and propagation matrices of modulus at most one, thanks to the $H_j^\pm$ values at the exponents in the outgoing expansions \eqref{eq:representation_up_layer}--\eqref{eq:representation_down_layer}, avoids the exponential growth of evanescent components that makes the transfer-matrix form unstable. 

\subsubsection{Energy flux identity} \label{s:flux_identity}
Let $u_j\in H^1_{\alpha_0}(\Omega_j)$ be a solution of the BVP \eqref{eq:bvp_single_layer} on the layer $\Omega_j$ for an incoming field $u_{\inc,j}=(u_{\inc,j}^\downarrow,u_{\inc,j}^\uparrow)$, and let  $\varepsilon$ be real-valued.
From the Green identity applied in each region with constant $\varepsilon$ and the regularity $\Delta u_j\in L^2(\Omega_j)$ (which prevents the presence of jumps in Dirichlet and Neumann traces of $u$ along the discontinuities of $\varepsilon$), we derive the energy flux identity \cite[Ch. 2]{ColtonKress2019}: 
\begin{equation} \label{eq:flux}
\int_{\partial\Omega_j} \Im \left\{\dn u_j \, \overline{u_j} \right\} \rd s = 0.
\end{equation}
The integrand vanishes on $\Gamma_D$, and the terms on the left and right boundaries cancel one another by quasi-periodicity, so that \eqref{eq:flux} reads
$\int_{\Gjp} \Im \{\dn u_j \, \overline{u_j}  \}\rd s+\int_{\Gjm} \Im \{\dn u_j \, \overline {u_j} \}\rd s=0$.
The expansion
$u_j =\sum_{n\in\IZ}  ( a_{j,n}^\downarrow \Psidown{j,n} + b_{j,n}^\uparrow \Phiup{j,n} )$
\eqref{eq:basis_expansion}
in $\Gjp \cup \Strip_j^+$ gives
\begin{align*}
u_j(x_1,H_j^+) = \sum_{n\in\IZ} \left( a_{j,n}^\downarrow + b_{j,n}^\uparrow \right) \re^{\ri \alpha_n x_1}, \qquad
\dxt u_j(x_1,H_j^+) = \sum_{n\in\IZ} \ri\beta_{j,n}^+ \left( b_{j,n}^\uparrow - a_{j,n}^\downarrow \right)\re^{\ri \alpha_n x_1}.
\end{align*}
The orthogonality $\int_{0}^{L} \re^{\ri (\alpha_n - \alpha_m)x_1}\rd x_1 = L\delta_{nm}$ gives
\begin{align*}
\int_{\Gjp} \Im \left\{\dn u_j \, \overline{u_j}  \right\}\rd s
&= \int_0^L \Im \big\{ \dxt u_j(x_1,H_j^+) \, \overline{u_j} (x_1,H_j^+) \big\}\rd x_1\\
&= L \sum_{n\in\IZ} \Im \left\{ \ri\beta_{j,n}^+ \left( b_{j,n}^\uparrow - a_{j,n}^\downarrow \right) \left( \overline{b_{j,n}^\uparrow} + \overline{a_{j,n}^\downarrow} \right) \right\}.
\end{align*}
\footnote{Using that $\varepsilon_j$ is assumed to be real, so $\beta_{j,n}^+$ is either real or purely imaginary by \eqref{eq:BetaEps}, we can further expand to
\begin{align*}
\int_{\Gjp} \Im \left\{\dn u_j \, \overline{u_j} \right\}\rd s
&=L \sum_{n\in\IZ, \;\alpha_n^2\le k^2\varepsilon_j^+}\beta_{j,n}^+
\big(| b_{j,n}^\uparrow|^2 - |a_{j,n}^\downarrow|^2\big)
-2 L\sum_{n\in\IZ, \;\alpha_n^2> k^2\varepsilon_j^+}|\beta_{j,n}^+|
\Im\{b_{j,n}^\uparrow \overline{a_{j,n}^\downarrow}\}.
\end{align*}}
Analogously, on the lower face
\begin{align*}
\int_{\Gjm} \Im \left\{\dn {u_j}  \, \overline{u_j}  \right\}\rd s
&= \int_0^L \Im \big\{ -\dxt u_j(x_1,H_j^-) \, \overline{u_j} (x_1,H_j^-) \big\}\rd x_1\\
&= L \sum_{n\in\IZ} \Im \left\{ \ri\beta_{j,n}^- \left( b_{j,n}^\downarrow - a_{j,n}^\uparrow \right) \left( \overline{b_{j,n}^\downarrow} + \overline{a_{j,n}^\uparrow} \right) \right\}.
\end{align*}
Plugging these in the flux identity \eqref{eq:flux}, we have
\begin{equation} \label{eq:flux_check}
\sum_{n\in\IZ} \Re
\left\{\beta_{j,n}^+ \left( b_{j,n}^\uparrow - a_{j,n}^\downarrow \right) \left( \overline{b_{j,n}^\uparrow} + \overline{a_{j,n}^\downarrow} \right)
+\beta_{j,n}^- \left( b_{j,n}^\downarrow - a_{j,n}^\uparrow \right) \left( \overline{b_{j,n}^\downarrow} + \overline{a_{j,n}^\uparrow}\right) \right\}
=0.
\end{equation}
This equality, together with the relation \eqref{eq:tmat} between $\mathbf{a}_j^{\uparrow\downarrow}$ and $\mathbf{b}_j^{\uparrow\downarrow}$, can be used to check the accuracy of a T-matrix implementation. 
We use this identity in \S\ref{s:ex2} to monitor the numerical consistency of the approximated T-matrix.

\section{DtN-TDG method for T-matrix approximation}\label{s:TDG}

\subsection{DtN-TDG method for scattering problems on a single layer}

We now present the DtN-TDG method, which is used to approximate the T-matrix $\Tmat_j^{M}$ \eqref{eq:tmat} of the BVP \eqref{eq:bvp_single_layer} posed on a single layer $\Omega_j$.
This method has been applied to such a periodic scattering problem with a plane wave incident field in \cite{Monforte2025}.
Here, we slightly extend the formulation in \cite{Monforte2025} to allow for the Fourier modes $\Psidown{j,n}$ and $\Psiup{j,n}$ as incident fields $u_{\inc,j}^\downarrow$ and~$u_{\inc,j}^\uparrow$.

Let $\Mh = \{ K \}$ be a convex-polygonal finite element partition of $\Omega_j$, such that the permittivity $\varepsilon$ is constant in each element.
We write $\bn_K$ for the outward-pointing unit normal vector on $\partial K$, $h_K$ for the diameter of $K$, and $h = \max_{K\in\Mh}h_K$ for the mesh width of $\Mh$.
We denote by $\nabla_h$ the elementwise application of $\nabla$ and write $\dn = \bn \cdot \nabla_h$ and $\dnK = \bn_K \cdot \nabla_h$ for the normal derivatives on $\partial\Omega_j$ and $\partial K$ respectively.

We require the mesh to be periodic: for each element $K\in\Mh$ with a face $F\subset\Gleft$, denoting its endpoints $(0,x_2^-)$ and $(0,x_2^+)$, there is an element $K'\in\Mh$ with a face $F'\subset\Gright$ with endpoints $(L,x_2^-)$ and $(L,x_2^+)$. 
Then $F$ and $F'$ are {identified and} treated as a single internal face, and $K$ and $K'$ are treated as adjacent elements. 

We denote by $\Fh = \bigcup_{K \in \Mh } \partial K $ and $\FhI = \Fh \setminus(\GD\cup\Gamma_{H^+}\cup\Gamma_{H^-})$ the skeleton of the mesh and its inner part.
Note that the faces in $\Gleft\cup\Gright$ are part of $\FhI$.
Given two adjacent elements $K_1, K_2\in \Mh$, an elementwise-smooth function $v$ and vector field $\pmb\tau$ on $\Mh$, we introduce on $\partial K_1 \cap \partial K_2$ the averages and the normal jumps:
\begin{align} \label{eq:Jumps}\begin{aligned}
\dlgraffa v\drgraffa & := \frac{1}{2}\left( v_{|K_1} + v_{|K_2} \right), & 
\dlgraffa\pmb\tau\drgraffa & := \frac{1}{2}\left( \pmb\tau_{|K_1} + \pmb\tau_{|K_2} \right), \\ 
\llbracket v \rrbracket_N & := v_{|K_1}\bn_{K_1} + v_{|K_2}\bn_{K_2},\qquad & 
\llbracket \pmb\tau \rrbracket_N & := 
\pmb\tau_{|K_1}\cdot \bn_{K_1} +  \pmb\tau_{|K_2}\cdot \bn_{K_2}.
\end{aligned}\end{align}
These face quantities must take into account quasi-periodicity:
if $K_1$ has two vertices with coordinates $(0,x_2^-)$ and $(0,x_2^+)$, and
$K_2$ has two vertices with coordinates $(L,x_2^-)$ and $(L,x_2^+)$, when writing the average and jump formulae \eqref{eq:Jumps} on their common face $v_{|K_1}$ and $\pmb\tau_{|K_1}$ have to be replaced by $v_{|K_1}\re^{\ri\alpha_0L}$ and $\pmb\tau_{|K_1}\re^{\ri\alpha_0L}$.

We then define the Trefftz space $T(\Mh)$, which contains both the BVP solution and the TDG basis functions:
\begin{align} \label{eq:BrokenSob}
T(\mathcal{M}_h) :=\big\{ & v \in L^2(\Omega_j) : \; v_{| K}\in H^1(K) \hspace{0.2cm} \forall K \in \mathcal{M}_h, \\
& \Delta v + \kappa^2v = 0 \hspace{0.1cm} \text{in} \hspace{0.1cm} K \hspace{0.1cm} \text{and} \hspace{0.1cm} \partial_{\mathbf{n}_K}v \in L^2(\partial K) \hspace{0.2cm}\forall K \in \mathcal{M}_h 
\big\}. \nonumber
\end{align}
The discrete Trefftz space $V_p(\Mh)$ is a finite-dimensional subspace of $T(\Mh)$ and can be represented as $V_p(\Mh) = \bigoplus_{K \in \Mh}V_{p_K}(K)$, where $V_{p_K}(K)$ is a $p_K$-dimensional subspace of $T(\Mh)$ of functions supported in $K$.

We use the TDG formulation of \cite{Hiptmair2011, Monforte2025}, which relies on the definition of ``numerical fluxes'' on the mesh skeleton.
Here, we omit the derivation of the TDG scheme and the definition of the numerical fluxes on $\FhI$ and $\Gamma_D$, which are identical to \cite{Hiptmair2011} and subsequent works.
On the horizontal boundaries $\Gjpm$, following \cite{Kapita2018,Monforte2025}, we define the numerical fluxes as
\begin{align*}
& \begin{cases}
\hat{u}_h = u_h -\ri\kappa^{-1}\pd\left(\nabla_h  u_h \cdot \bn - \DtN_{\varepsilon_j^+,M}^+ u_h + 2\ri\beta_{j,n}^+ u_{\inc,j}^\downarrow\right), 
\hspace{30mm} \text{on }\Gjp,\\
\widehat{\ri \kappa \sigma}_h = -  \DtN_{\varepsilon_j^+,M}^+ u_h \bn + 2\ri\beta_{j,n}^+ u_{\inc,j}^\downarrow \bn 
\\ \hspace{20mm}
+\ri\kappa^{-1}\pd \DtN^{+,*}_{\varepsilon_j^+,M} \left( \nabla_h  u_h -\DtN_{\varepsilon_j^+,M}^+ u_h \bn + 2\ri\beta_{j,n}^+ u_{\inc,j}^\downarrow \bn \right),
\end{cases}
\\
& \begin{cases}
\hat{u}_h = u_h -\ri\kappa^{-1}\pd\left( \nabla_h  u_h \cdot \bn - \DtN_{\varepsilon_j^-,M}^- u_h + 2\ri\beta_{j,n}^- u_{\inc,j}^\uparrow \right),
\hspace{30mm} \text{on }\Gjm,
\\
\widehat{\ri \kappa \sigma}_h = - \DtN_{\varepsilon_j^-,M}^- u_h \bn + 2\ri\beta_{j,n}^- u_{\inc,j}^\uparrow \bn
\\ \hspace{20mm}
+ \ri\kappa^{-1}\pd \DtN^{-,*}_{\varepsilon_j^-,M} \left( \nabla_h  u_h -\DtN_{\varepsilon_j^-,M}^- u_h \bn + 2\ri\beta_{j,n}^- u_{\inc,j}^\uparrow \bn\right),
\end{cases}
\end{align*}
where $\pd \in L^\infty(\Gjp\cup\Gjm)$ is a positive flux coefficient, and 
$\DtN^{\pm,*}_{\varepsilon_j^\pm,M}$ is the $L^2(\Gjpm)$-adjoint of $\DtN^\pm_{\varepsilon^\pm,M}$:
\begin{equation*}
\int_{\Gjpm} \DtN^{\pm,*}_{\varepsilon^\pm,M} v \, \cw \rd s
= \int_{\Gjpm} v \, \overline{\DtN^\pm_{\varepsilon^\pm,M} w} \rd s
\qquad \forall v,w\in H^{1/2}_{\alpha_0}({\Gjpm}), \qquad M \in \IN .
\end{equation*}

The terms appearing in the numerical fluxes on $\Gjpm$ are consistent, i.e.\ $\hat u_h=u_h$ and $\bn\cdot\widehat{\ri \kappa \sigma}_h=\partial_\bn u_h$ if $u_h$ satisfies the boundary conditions in \eqref{eq:bvp_single_layer}, because, for the Fourier modes defined in \eqref{eq:psi_phi},
\begin{align*}
& \dn \Psidown{j,n} = -\ri \beta_{j,n}^+ \Psidown{j,n}, \qquad 
\DtN_{\varepsilon_j^+}^+ \Psidown{j,n} 
= \ri \beta_{j,n}^+ \Psidown{j,n}, \qquad \, \text{on} \,\, \Gjp, \\
& \dn \Psiup{j,n} = -\ri \beta_{j,n}^- \Psiup{j,n}, \qquad \DtN_{\varepsilon_j^-}^- \Psiup{j,n} 
= \ri \beta_{j,n}^- \Psiup{j,n}, \qquad \, \text{on} \,\, \Gjm.
\end{align*}

The DtN-TDG formulation reads as follows: Find $u_h^M \in V_p(\Mh)$ such that
\begin{equation}\label{eq:TruncDG}
\cA_h^M(u_h^M, v_h) =\ell_h^M(v_h) 
\qquad \forall v_h \in V_p(\Mh), \qquad M\in \IN, 
\end{equation}
where
\begin{align*}\nonumber
&\cA_h^{ M}(u, v) :=\\
& \int_{\FhI} \left( \dlgraffa u\drgraffa \llbracket \overline{\nabla_h v} \rrbracket_N  
-  \dlgraffa\nabla_h u\drgraffa \cdot \llbracket \cv \rrbracket_N 
- \ri\xi \,  {\pa} \, \llbracket  u\rrbracket_N \cdot \llbracket \cv \rrbracket_N    
- \ri\xi^{-1} {\pb} \, \llbracket \nabla_h u\rrbracket_N  \, \llbracket \overline{\nabla_h v} \rrbracket_N \right) \rd s \\ 
& + \int_{\Gjp} \biggl( u \, \overline{\dn v }  - \DtN_{\varepsilon_j^+,M}^+u \, \overline{v}
-{\pd} \, \ri\kappa^{-1} \left(\dn u -\DtN_{\varepsilon_j^+,M}^+ u\right) \, \overline{\left(\dn v -\DtN_{\varepsilon_j^+,M}^+ v\right)} \biggl) \rd s \\ \nonumber
& + \int_{\Gjm} \biggl(  u \, \overline{\dn v}  - \DtN_{\varepsilon_j^-,M}^-u \, \overline{v}
-{\pd} \, \ri\kappa^{-1} \left(\dn u -\DtN_{\varepsilon_j^-,M}^- u\right) \, \overline{\left(\dn v -\DtN_{\varepsilon_j^-,M}^- v\right)} \biggl) \rd s \\ \nonumber
& + \int_{\GD\cap\partial\Omega_j} \left( -\dn u \, \cv - \ri\kappa \,{\pa} \, u\, \cv \right) \rd s,
\end{align*}
where  $\pa\in L^\infty(\FhI\cup \Gamma_D)$ and $\pb\in L^\infty(\FhI)$ are positive flux coefficients, and
\begin{equation*}
\xi: = \frac{\Re (\kappa_{|K_1}) + \Re (\kappa_{|K_2})}2
\qquad \text{on the internal face }F=\partial K_1\cap\partial K_2,\quad K_1,K_2\in\Mh.
\end{equation*}
The functional $\ell_h^M$ is 
\begin{align*} 
\ell_h^M(v) :=  & \int_{\Gjp} - 2\ri\beta_{j,n}^+ \, u_{\inc,j}^\downarrow \left( \cv - {\pd} \, \ri\kappa^{-1} \, \overline{ (\partial_\bn v -\DtN_{\varepsilon_j^+,M}^+ v )} \, \right) \rd s \\ \nonumber
& + \int_{\Gjm} - 2\ri\beta_{j,n}^- \, u_{\inc,j}^\uparrow\left( \cv - {\pd} \, \ri\kappa^{-1} \, \overline{ (\partial_\bn v -\DtN_{\varepsilon_j^-,M}^- v )} \, \right) \rd s.
\end{align*}  

Since the faces on $\Gamma_{\rm left}\cup\Gamma_{\rm right}$ are treated as internal faces, with the averages and jumps modified by the factor $\re^{i\alpha_0 L}$ as described in the paragraph after \eqref{eq:Jumps}, the formulation \eqref{eq:TruncDG} enforces the quasi-periodicity of the solution.
The DtN-TDG formulation \eqref{eq:TruncDG} differ from \cite[eq.~(45)]{Monforte2025} only in the right-hand side $\ell_h^M$.
The consistency of the numerical fluxes implies the consistency of the method.
The coercivity and the quasi-optimality of the DtN-TDG method proved in \cite[Prop.~4.2]{Monforte2025} apply to \eqref{eq:TruncDG} under the same non-trapping assumptions and for $M\ge M_\star$ in \eqref{eq:MLarge}.

The most common TDG schemes (see \cite{Hiptmair2011,Hiptmair2016}), including our implementation, 
use plane waves as basis functions of the TDG method.
For a mesh element $K \in \Mh$, we denote by $V_p(K)$ the plane wave space on $K$ spanned by $p$ plane waves, $p \in \mathbb{N}$:
\begin{equation*}
V_p(K) = \bigg\{ \; v \in L^2(K) : \; v(\bx) = \sum_{\ell=1}^p \eta_\ell
\exp\{\ri\kappa\bd_\ell \cdot \bx \}, \hspace{0.3cm} \eta_\ell \in \mathbb{C} \; \bigg\},
\end{equation*}
where $\{ \bd_\ell \}_{\ell=1}^p \subset \mathbb{R}^2$, with $|\bd_\ell|=1$, are different propagation directions.
Isotropic discrete spaces are achieved by distributing the directions uniformly on the unit circle, i.e.\ choosing
$\bd_\ell = ( \cos\frac{2\pi\ell}p, \sin\frac{2\pi\ell}p)$, $\ell=1,\ldots,p$. 
For simplicity, we select the same number $p$ of directions in every element $K\in\Mh$.
The value of $\kappa$ depends on the region where the element $K$ is located;
recall that we consider meshes such that 
$\kappa$ is constant inside each element.
The corresponding global discrete space $V_p(\Mh)$ is
\begin{equation*}
V_p(\Mh):=\bigoplus_{K\in\Mh}V_p(K)
=\big\{ \; v \in L^2(\Omega) : v_{|K} \in V_p(K), \;\forall K \in \Mh \; \big\}.
\end{equation*}
An advantage of plane-wave bases is that all entries of the stiffness matrix and right-hand side vector admit explicit analytical expressions, eliminating numerical quadrature errors and reducing assembly cost. The derivation of these entries is described in detail in \cite[\S4.3]{Monforte2025}.

Plane-wave discrete spaces allow to achieve highly accurate approximations.
In particular, for isotropic plane wave spaces, convergence for $p\to\infty$ is exponentially fast when the target solution admits an analytic extension on a neighbourhood of the mesh element, see \cite{Hiptmair2011,Hiptmair2016}.
On the other hand, the plane-wave approximation of certain Helmholtz solutions requires very large coefficients in their linear combinations, leading to numerical cancellation and strong instabilities.
These instabilities might be remedied by using discrete spaces containing evanescent plane waves, see \cite{PaHuMo2023,MoGaPa2026}.

\begin{rem}[Transverse-magnetic (TM) waves]\label{rem:TMmodes}
We have considered the case of TE-polarized waves.
TM-polarized waves can be treated similarly.
In this case, the PEC obstacle would be represented by a Neumann boundary condition ($\partial_\bn u=0$ on $\partial D$) and the discontinuous coefficient would appear in the second-order term of the PDE operator as $\nabla\cdot(\frac1{\varepsilon}\nabla u)+k^2 u=0$, see \cite{BonnetBD1994}.
Both differences are easily dealt with the TDG discretization, see \cite{Huttunen2002}.
However, the BVP solution $u$ might have stronger singularities at material corner interfaces.
Non-trapping conditions for TM-polarized waves can be found e.g.\ in \cite{Civiletti2020}.
\end{rem}

\subsection{From the TDG solution to the T-matrix}
The T-matrix $\Tmat_j^{ M}$ \eqref{eq:tmat} maps the Fourier coefficient vector $(\aup_j,\adown_j)$ of any incident field to the coefficients vector $(\bup_j,\bdown_j)$ of the corresponding scattered field.
To approximate $\Tmat_j^{ M}$ column-by-column, we follow two steps.

First, we solve $4M+2$ times the BVP \eqref{eq:bvp_single_layer} using the DtN-TDG scheme \eqref{eq:TruncDG} with incident fields 
$(u_{\inc,j}^\uparrow,u_{\inc,j}^\downarrow)=(0,\Psidown{n})$ and 
$(u_{\inc,j}^\uparrow,u_{\inc,j}^\downarrow)=(\Psiup{n},0)$, for $n\in\{-M,\ldots,M\}$.
These $4M+2$ linear systems only differ from one another in the right-hand side vector: the TDG stiffness matrix needs to be assembled once. 
Furthermore, if an LU (or other) decomposition is used to solve the resulting linear system, then this needs to be computed only once.

Then, for each incoming field, the Fourier coefficients on $\Gjpm$ of the TDG solution have to be computed.
Let the DtN-TDG solution $u_h^M$ of \eqref{eq:TruncDG} be $u_h^M|_K=\sum_{\ell=1}^p\eta_\ell^K\,\varphi_\ell^K$ on each element $K\in\Mh$, where $\varphi_\ell^K(\bx)=\re^{\ri\kappa\bd_\ell\cdot\bx}$ are the TDG basis functions supported in the element $K$ and $\eta_\ell^K$ the relevant entries of the Galerkin solution vector.
Let $K\in\Mh$ be a mesh element adjacent to $\Gjp$ with $\partial K\cap\Gjp=[a,b]\times\{H_j^+\}$.
Then, thanks to formula \eqref{eq:FourierCoef}, the Fourier coefficients of $\varphi_\ell^K$ on $\Gjp$ are 
$$
\varphi_{\ell,n}^K(H_j^+)=\frac1L\re^{\ri\kappa(\bd_\ell)_2H_j^+}\int_a^b\re^{\ri(\kappa(\bd_\ell)_1-\alpha_n)x_1}\rd x_1,\qquad n\in\{ -M, \ldots, M\},
$$
which can be computed analytically as in \cite[\S4.3]{Monforte2025}.
Thus, the Fourier coefficients of $u_h^M$ are
$$
u_{h,n}^M(H_j^+)=\sum_{K\text{ adjacent to }\Gjp} \sum_{\ell=1}^p\eta_\ell^K \varphi_{\ell,n}^K(H_j^+),
\qquad n\in\{ -M, \ldots, M\}.
$$
Moreover, the coefficients $\varphi_{\ell,n}^K (H_j^\pm)$ are already required in the assembly of the DtN-TDG stiffness matrix, so there is no need to compute them again.

Once derived, the Fourier coefficients $u_{h,n}^M(H_j^\pm)$ are used as entries of the T-matrix $\Tmat_j^{M}$, after subtracting the expansion of the incoming field $\Psi_{j,n}^\downarrow$ or $\Psi_{j,n}^\uparrow$, which has only one non-zero Fourier component.
This shows that, given the DtN-TDG solution vectors, all the entries of the approximated and truncated T-matrix 
are computed analytically with no need for quadrature formulas.

\begin{rem}[Approximation of the solution inside the layers]\label{rem:PostProc}
The coupled linear system \eqref{eq:coupling_system} returns the Fourier coefficients $\bup_j,\bdown_j$ of the outgoing components of the approximate solution $u$ of the global BVP \eqref{eq:BVP} on the layer faces $\Gjpm$.
The expansions \eqref{eq:Representation} and \eqref{eq:ReprDown} immediately allow to compute the approximation of $u$ in the complement of the layers, i.e.\ in the semi-infinite strips $\Strip^\pm$ and in the gaps $(0,L)\times[H_j^+,H_{j-1}^-]$.
To evaluate the solution in each layer $\Omega_j$, in a post-processing stage we solve a local BVP on $\Omega_j$ itself using the TDG scheme, with the same mesh and discrete space already used in the T-matrix computations.
To ensure the stability of these local BVPs, on $\Gjpm$ we impose the impedance condition $\nabla u \cdot \bn - \ri \kappa u = g_R$, where the datum $g_R$ is derived from the 
expansion \eqref{eq:basis_expansion} of $u$ on $\Gjpm$, which is computable from the solution vector of the T-matrix system \eqref{eq:coupling_system}.
The derivatives of the basis functions $\Phidown{n}$ and $\Phiup{n}$ enter the impedance boundary datum $g_R$ and are exponential functions, so all matrix and vector entries are computed with a closed formula. 
The TDG numerical fluxes on impedance boundaries are defined in \cite[\S2]{Hiptmair2011}.
\end{rem}

\begin{rem}[Computational complexity analysis]\label{rem:OpCount}
We compare the complexity (in terms of asymptotic operation count) of the DtN-TDG scheme with and without the T-matrix approach, in a simple setting.
Let $N$ be the number of layers, each partitioned in a mesh of $E\approx E_1E_2$ polygonal elements, in such a way that there are roughly $E_1$ elements in the $x_1$ direction and $E_2$ elements in the $x_2$ direction (e.g.\ a $E_1\times E_2$ tensor-product mesh, possibly with some perturbation), and set $E_0:=\min\{E_1,E_2\}$.
Let $p$ be the dimension of the local Trefftz space on each mesh element, and $4M+2$ the dimension of each local T-matrix.
Assume that the layers are adjacent to one another, i.e.\ $H_j^-=H_{j+1}^+$.

Assume that $M\le\delta pE_1$ for some $0<\delta\le1$. 
The ratio $\frac M{pE_1}$ measures how much the retained Fourier expansion of the transmitted field is compressed compared to the representation in the TDG basis.

The monolithic TDG method has $O(pEN)$ degrees of freedom (DOFs), of which $O(pE_1)$ adjacent to $\Gamma_{H^\pm}$.
The number of the TDG matrix non-zeros is $O(p^2 EN)$ since each element is adjacent to a small number of other elements.
The TDG matrix assembly cost is $O(p^2EN + p^2 E_1 M)$, where the first term accounts for the bulk of the matrix and the second for the entries involving the truncated DtN map.
The monolithic system, if appropriately ordered, has bandwidth $O(pE_1)$, so its solution by a banded direct solver has cost $O((pE_1)^2 (pEN))$ \cite[\S4.3.1]{Golub2013}.
So the total cost is $O(p^2EN + p^2 E_1 M + p^3E_1^2EN)=O(p^3E_1^2EN)$.

Consider now the T-matrix approach.
For each layer: 
the TDG matrix has size $pE$ and has bandwidth $O(pE_0)$ if the elements are appropriately ordered;
its assembly has cost $O(p^2E+p^2E_1M)$;
its LU factorization $O(p^3E_0^2E)$;
the solution of $4M+2$ linear systems by forward and backward banded substitution $O(p^2E_0EM)$  \cite[\S4.3.3]{Golub2013};
the T-matrix assembly $O(pE_1M^2)$ (recalling that the Fourier coefficient of the plane wave basis functions are computed during the TDG matrix assembly);
its LU factorization $O(M^3)$;
the solution of $M$ dense linear systems by forward and backward substitution $O(M\cdot M^2)$.
The solution of the global block-tridiagonal T-matrix linear system costs $O(NM^3)$.

Thus the T-matrix technique costs 
$O(p^2E+p^2E_1M+p^3E_0^2E+p^2E_0EM+pE_1M^2+M^3+  M^3N)
=O(p^3E_0^2E+ \delta p^3E_1^2E+\delta^3p^3E_1^3N)$
if only a small finite number of different layers are present.
If $\delta\approx1$, i.e.\ $M\approx pE_1$ and the Fourier modes retained in the layer coupling are proportional to the TDG DOFs along the interface, then the the cost is $O(p^3E_1^2(E+E_1N))$ and is cheaper than the monolithic approach by a factor $O(\min\{E_2,N\})$.
If only fewer Fourier modes are retained ($\delta\ll1$, e.g.\ $E_1\le E_2$ and $\delta\le(E_2/N)^{1/3}$), then the T-matrix approach is cheaper by a factor $O(\max\{1,\frac{E_1}{E_2}\}^2N)\ge O(N)$.

Instead, if $O(N)$ layers are different from one another, $O(N)$ different TDG and $\mathbb T_j^M$ matrices need to be assembled and factorized.
Then, the cost of the T-matrix technique is
$O([p^2E+p^2E_1M+p^3E_0^2E+p^2E_0EM+pE_1M^2+M^3]N)
=O(p^3[E_0^2E+ \delta E_1^2E]N)$, which is cheaper than the monolithic approach by a factor $\min\{1/\delta,\max\{1,\frac{E_1}{E_2}\}^2\}$.
In this case, the T-matrix implementation is convenient, asymptotically for large parameter values, if the number of Fourier mode is moderate ($\delta\ll1$) \emph{and} the layers are "shallow", i.e.\ more extended in the $x_1$ than the $x_2$ direction ($E_1\gg E_2$).

The T-matrix enjoys a further considerable advantage if gaps are present between the layers and/or if parallel implementations are used; see the numerical example in \S\ref{s:ManyLayers}.
\end{rem}

\section{Numerical experiments}\label{s:Numerics}
We present some numerical experiments to assess the accuracy of the method proposed in the approximation of different BVPs, with  and without impenetrable obstacles.
In all these experiments, the T-matrix is approximated using the DtN-TDG method \eqref{eq:TruncDG} and the solution outside the layers is obtained solving \eqref{eq:coupling_system}.  
The DtN-TDG scheme has been implemented in MATLAB; all TDG linear systems are solved with the ``backslash'' direct solver, while the coupled system \eqref{eq:coupling_system} is solved with the block Band Gaussian Elimination algorithm (see \cite[\S4.3.1]{Golub2013}).
The triangulation on each layer $\Omega_j$ is made with the MATLAB PDE toolbox; to evaluate the $L^2(\Omega)$ and $H^1(\Omega)$ norms of the DtN-TDG error we use a Duffy quadrature rule on the mesh triangles, which is presented in \cite{Duffy1982} and implemented in \cite{Kubatko2024}. 
The numerical flux parameters used correspond to the ultra weak variational formulation (UWVF) of Cessenat and Després \cite{Cessenat1998}, i.e.\ $\pa=\pb=\pd=\frac{1}{2}$ (see \cite[\S2.2.2]{Hiptmair2016}).

\subsection{Two flat interfaces between homogeneous materials}\label{s:ex1}
\begin{figure}[htb]
\centering
\includegraphics[height=.31\textwidth]{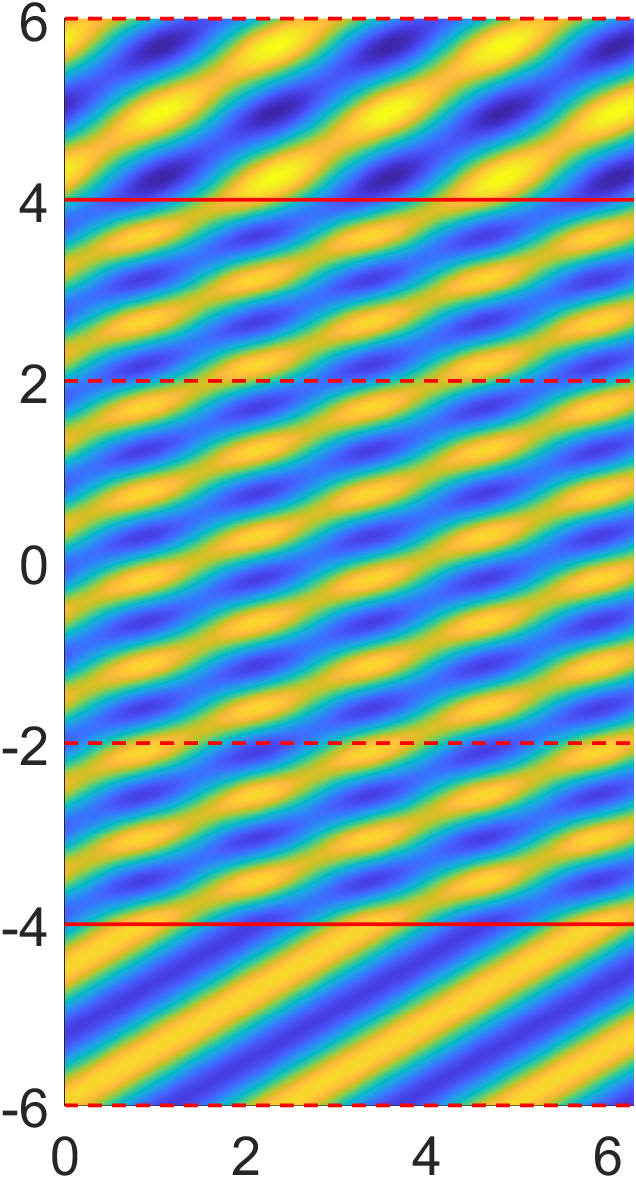}\hspace*{8pt}
\includegraphics[height=.31\textwidth]{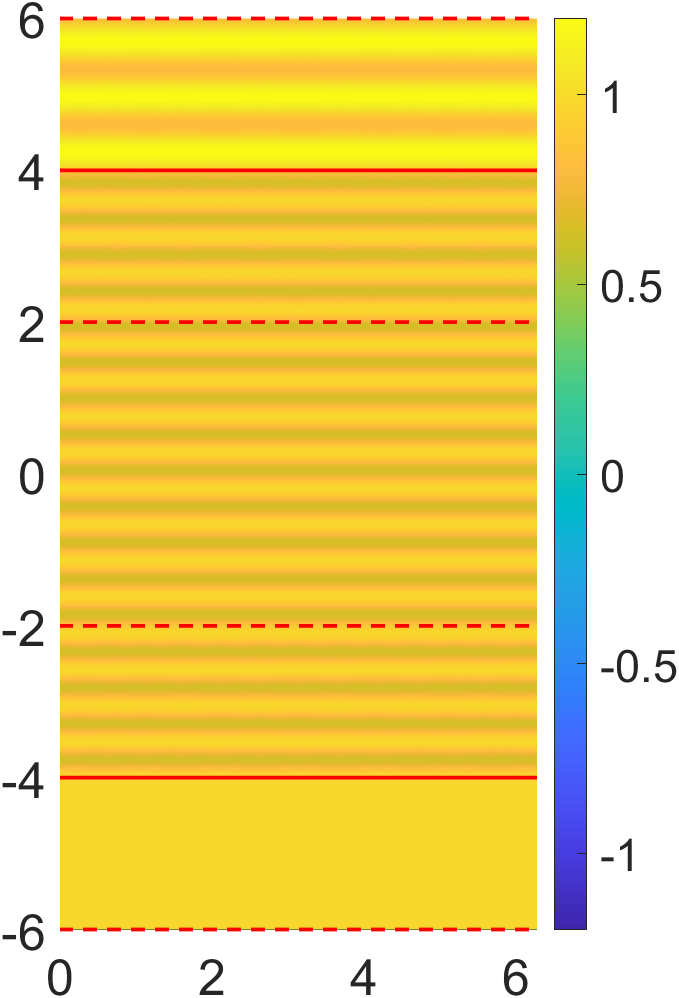}\hfill
\includegraphics[height=.31\textwidth]{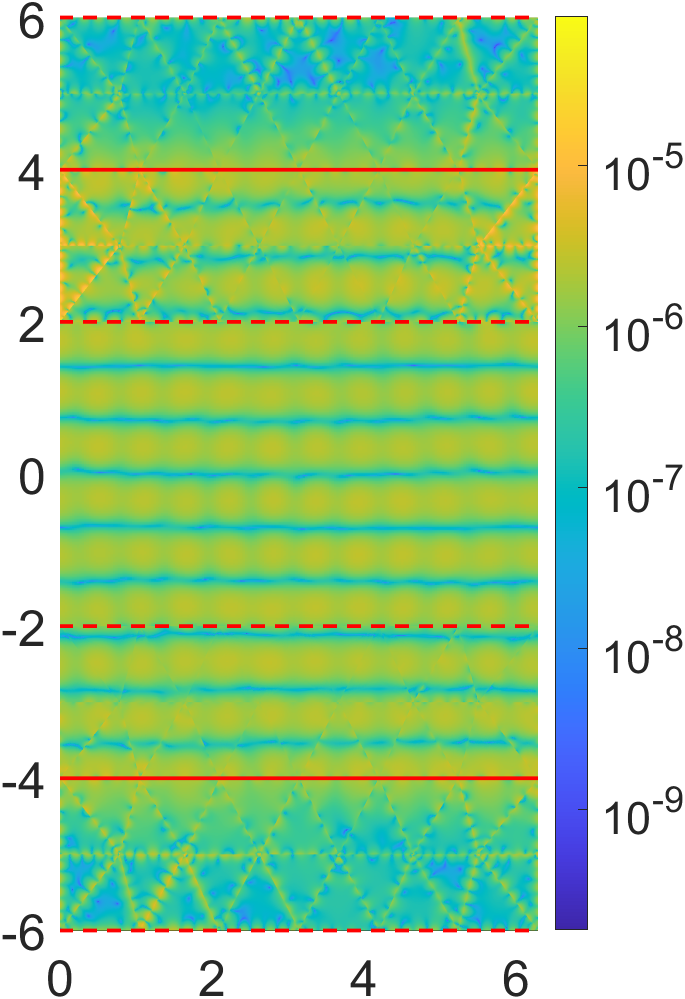}\hfill
\includegraphics[height=.27\textwidth]{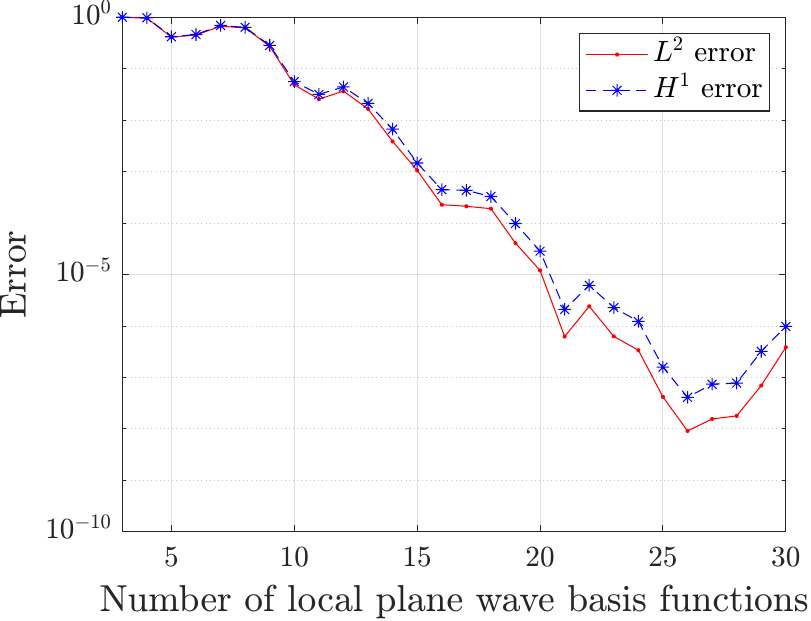}
\caption{Left to right: real part and absolute value of the numerical solution (with the same color bar) of the BVP in \S\ref{s:ex1}; 
absolute value of the pointwise error (in logarithmic color scale) for $h=1$, $M=18$ and $p=30$; convergence of the $L^2(\Omega)$ and $H^1(\Omega)$ relative error norms for $p \in \{3,\ldots,30\}$.
Continuous red lines indicate the interfaces between different materials, while dashed lines indicate layers boundaries.}
\label{fig:twolayers_eps2}
\end{figure}

We consider two layers $\Omega_1 = (0,2\pi)\times(H_1^-,H_1^+)$, $\Omega_2=(0,2\pi)\times(H_2^-,H_2^+)$, with $H_1^+ > H_1^- > 0$ and $H_2^+ = -H_1^-$, $H_2^- = - H_1^+$. $\Omega_1$ and $\Omega_2$ are both composed of two different materials, separated by the horizontal line $\{x_2=\frac12(H_j^+ + H_j^-)\}$, such that $\varepsilon_1^+ = \varepsilon_2^- = 1$ and $\varepsilon_1^-=\varepsilon_2^+=\varepsilon_{in}>1$. In this case, the total field in $\Omega = \Omega_1 \cup \Omega_2 \cup (0,2\pi)\times[H_2^+,H_1^-]$ can be computed analytically and we can compute the error of the numerical solution against the exact one.
Defining $d=\frac12(H_1^+ + H_1^-)$ and $\gamma := \sqrt{\varepsilon_{in} - \cos^2 \theta}$, the total field in $\Omega$ is:
\begin{equation*}
u(x_1,x_2)  = \begin{cases} 
\re^{\ri k(x_1 \cos \theta + x_2 \sin \theta)} + \mathtt{R}\re^{\ri k(x_1 \cos \theta - x_2\sin\theta)} & x_2>d,\\
\mathtt{T_1}\re^{\ri k(x_1 \cos \theta - x_2 \gamma)} + \mathtt{T_2}\re^{\ri k(x_1 \cos \theta + x_2 \gamma)} & -d<x_2<d, \\
\mathtt{T_3}\re^{\ri k(x_1 \cos \theta + x_2 \sin \theta)} & x_2<-d. \end{cases}
\end{equation*}
By enforcing the continuity of $u$ and $\partial_{x_2} u$ at $x_2 = d$ and $x_2 = -d$, we obtain the following 4-dimensional linear system for the reflection and transmission coefficients $\mathtt{R}$, $\mathtt{T_1}$, $\mathtt{T_2}$, and $\mathtt{T_3}$:
\begin{equation*}
\begin{cases}
\re^{-\ri kd\sin\theta} \mathtt{R} - \re^{-\ri kd \gamma} \, \mathtt{T_1} - \re^{\ri kd \gamma} \, \mathtt{T_2}
= -\re^{\ri kd\sin \theta},\\
-\sin \theta\, \re^{-\ri kd\sin\theta} \mathtt{R} + \gamma \re^{-\ri kd \gamma} \, \mathtt{T_1} 
- \gamma \re^{\ri kd \gamma}\, \mathtt{T_2} = -\sin \theta \,\re^{\ri kd\sin\theta}, \\
\re^{\ri kd \gamma} \, \mathtt{T_1} + \re^{-\ri kd \gamma} \, \mathtt{T_2} - \re^{-\ri kd\sin\theta} \, \mathtt{T_3} = 0, \\
-\gamma \re^{\ri kd \gamma} \, \mathtt{T_1} + \gamma \re^{-\ri kd \gamma} \, \mathtt{T_2} - \sin \theta \, \re^{-\ri kd\sin\theta} \, \mathtt{T_3} = 0.
\end{cases}
\end{equation*}
We consider an incident plane wave with wavenumber $k=5$ and incident angle $\theta=-\pi/3$, and we choose $H_1^+=6$ and $H_1^-=2$, so that $d=4$.
We select $\varepsilon_{in}=2$ and a mesh with 52 triangles for each layer $\Omega_1$ and $\Omega_2$.
We study the $p$-convergence of the method, increasing the number $p$ of plane wave basis functions in every element of the mesh.
The results are displayed in Figure \ref{fig:twolayers_eps2}: as expected, we observe exponential convergence of the $L^2(\Omega)$ and $H^1(\Omega)$ errors, until the plane wave basis instability kicks in, \cite{PaHuMo2023,MoGaPa2026,Hiptmair2016}.

\subsection{Interfaces with corners and convergence in the DtN truncation parameter \texorpdfstring{$M$}M}\label{s:ex2}
In this section, we apply the DtN-TDG scheme to a problem involving corner singularities.
The exact solution is not available in closed form, so we compute the error by testing against a DtN-TDG solution computed with higher numbers of plane waves on the same mesh.

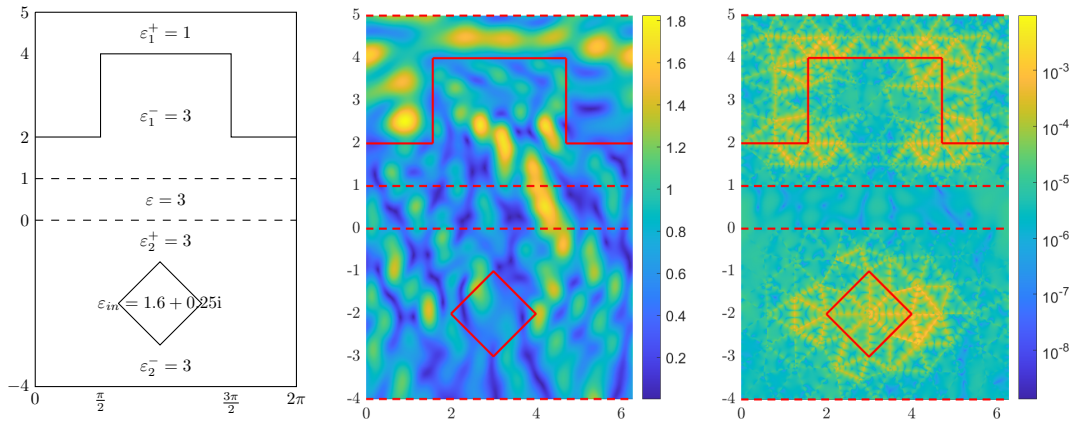
\begin{figure}[htbp]
\centering
\begin{tikzpicture}[scale=.55, every node/.style={scale=0.65}]
\draw (0,5) node[left]{$5$};
\draw (0,4) node[left]{$4$};
\draw (0,2) node[left]{$2$};
\draw (0,1) node[left]{$1$};
\draw (0,0) node[left]{$0$};
\draw (0,-4) node[left]{$-4$};
\draw (0,-4) node[below]{$0$};
\draw (1.57, -4) node[below]{$\frac{\pi}{2}$};
\draw (4.71, -4) node[below]{$\frac{3\pi}{2}$};
\draw (6.28, -4) node[below]{$2\pi$};
\node[] at (3.14,4.5) {$\varepsilon_1^+ = 1$};
\node[] at (3.14,2.5) {$\varepsilon_1^- = 3$};
\draw (0,5) -- (6.28,5) -- (6.28,-4) -- (0,-4) -- (0,5); 
\draw (0,2) -- (1.57,2) -- (1.57,4) -- (4.71,4) -- (4.71,2) -- (6.28,2);
\draw[dashed] (0,1) -- (6.28,1);
\draw[dashed] (0,0) -- (6.28,0);
\draw (2,-2) -- (3,-1) -- (4,-2) -- (3,-3) -- (2,-2);
\node[] at (3,-2) {$\varepsilon_{in} = 1.6+0.25\ri$};
\node[] at (3.14,-0.5) {$\varepsilon_2^+ = 3$};
\node[] at (3.14,-3.5) {$\varepsilon_2^- = 3$};
\node[] at (3.14,0.5) {$\varepsilon = 3$};
\end{tikzpicture}\hfill
\includegraphics[height=.385\textwidth]{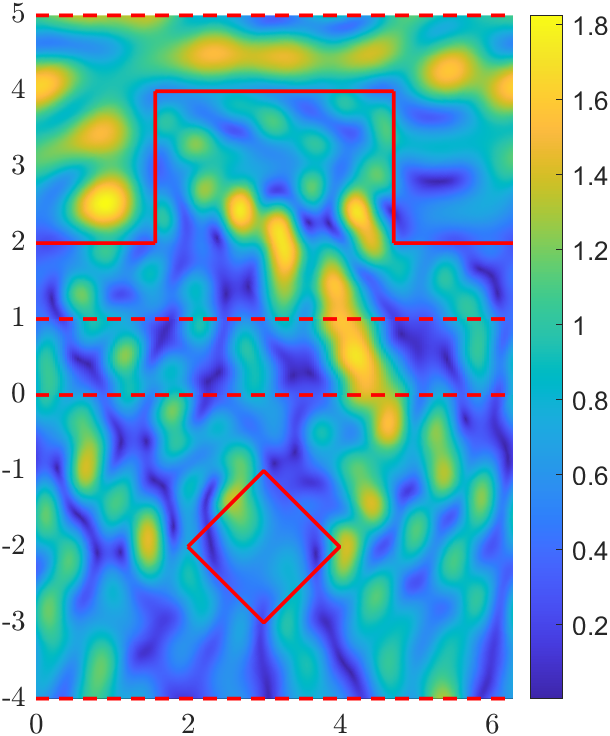}\hfill
\includegraphics[height=.385\textwidth]{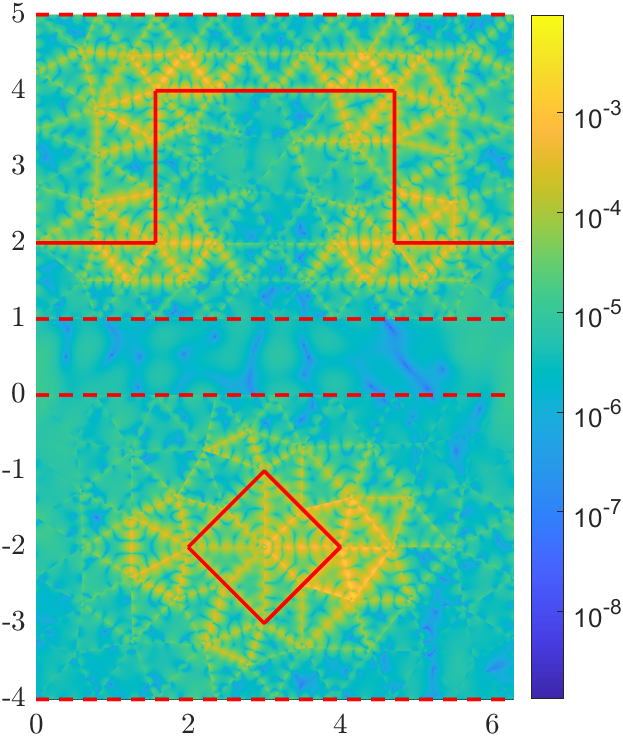}
\caption{The problem with corner singularities in \S\ref{s:ex2}.
Left: the geometry of the problem.
Center: the absolute value of the numerical solution.
Right: the absolute value of the TDG error for $p=19$ and $M=15$ against a refined numerical solution, in logarithmic color scale.}
\label{fig:corner_sing}
\end{figure}

We consider the domain $\Omega = (0, 2\pi) \times (-4, 5)$ containing two stacked layers $\Omega_1 = (0,2\pi) \times (1,5)$ and $\Omega_2 = (0,2\pi) \times (-4,0)$, as shown in Figure \ref{fig:corner_sing}. In the first layer the relative permittivity assumes a value $\varepsilon=\varepsilon_1^+ = 1$ in $\{x_2>4\}\cup\{x_2>2, |x_1-\pi|>\frac\pi2\}$ and $\varepsilon=\varepsilon_1^- = 3$ otherwise, while in the second layer the relative permittivity assumes the value $\varepsilon = \varepsilon_{in} = 1.6+0.25\ri$ inside the square of vertices $(2,-2)$, $(3,-1)$, $(4,-2)$ and $(3,-3)$, and $\varepsilon=\varepsilon_2^+=\varepsilon_2^-=3$ otherwise. In the region between the layers the relative permittivity assumes the value $\varepsilon=3$.
We select an incident plane wave with wavenumber $k = 4$ and incident angle $\theta=-\pi/3$.
The meshes on $\Omega_1$ and $\Omega_2$ contain 118 and 104 triangles, respectively.

Figure \ref{fig:corner_sing} shows the approximate solution and the pointwise error, using as reference a DtN-TDG solution with $p=20$, showing how the error concentrates on the boundaries of the elements adjacent to the material interface corners.
The scatterer profile is also displayed in the solution plot.

\begin{figure}[htbp]
\centering
\includegraphics[height=.35\textwidth]{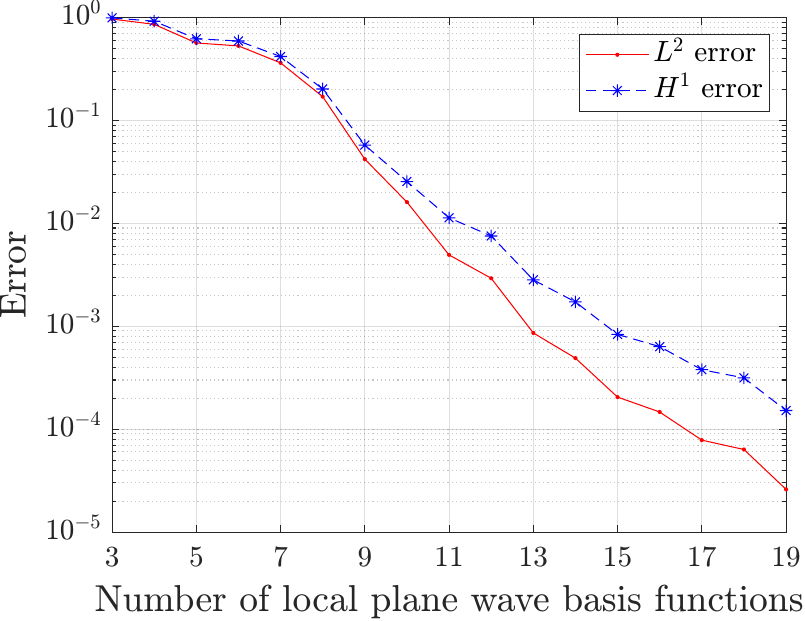}\hfill
\includegraphics[height=.35\textwidth]{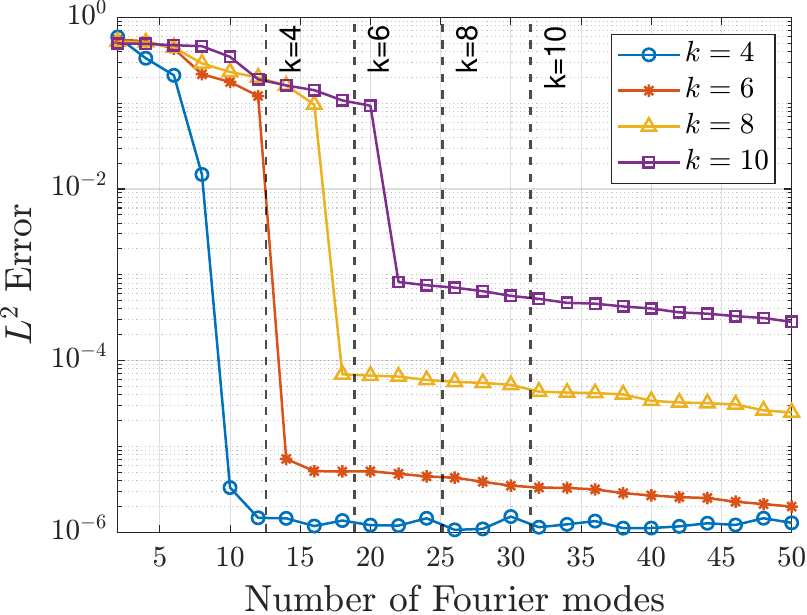}
\caption{Left: the $p$-convergence of the relative $L^2(\Omega)$ and $H^1(\Omega)$ errors with $M=15$ for the example in \S\ref{s:ex2}.
Right: relative $L^2(\Omega)$ error against the truncation order $M$ of the DtN operator, for $k= 4, 6, 8, 10$; all experiments use $p=20$ plane waves per element.
Dashed vertical lines correspond to the value $kL/2$.}
\label{fig:pConv_MConv}
\end{figure}

In the left plot of Figure~\ref{fig:pConv_MConv}, we show the decay in $p$ of the relative $L^2(\Omega)$ and $H^1(\Omega)$ errors on a fixed mesh, with $M=15$.
In the right plot, we show the error in dependence of the Fourier parameter $M$ for different values of the wavenumber $k$ and $p=20$; the same value of $p$ is used for the reference solution.
Recall that, in each experiment, the same value of $M$ is used to truncate the T-matrix blocks (\S\ref{s:TmatApprox}) and the DtN map in the TDG fluxes (\S\ref{s:TDG}).
The error decreases quickly in $M$ and then flattens, showing that the Fourier-truncation error is dominated by the TDG discretization.
As expected, larger values of $M$ are needed for larger wavenumbers $k$.
This is in complete agreement with the DtN-TDG results for the periodic scattering problem in \cite[Fig.~12]{Monforte2025} and for bounded scatterers \cite[Fig.~2]{Kapita2018}.
We observe that the relative error for $k=4$ in the $M$-convergence plot is lower than the corresponding relative error in the $p$-convergence plot.
This suggests that, once $M$ is higher than some threshold value, increasing
$p$ has a more significant impact on the numerical error than increasing the number of Fourier modes $M$.

\begin{rem}[Truncation order dependence on parameters] \label{rem:truncation}
Repeating the M-convergence experiment with the same parameters and a smaller period $L$, we observe that the order of truncation $M$ required to reach the error plateau is smaller, while changing the vertical width of the layers does not affect the results. Experimentally, we derive the following rule for the truncation order:
\begin{equation*}
M > \frac{k L}{2},
\end{equation*}
and we apply it in our numerical experiments.
\end{rem}

\begin{figure}[htbp]
\centering
\includegraphics[height=.37\textwidth]{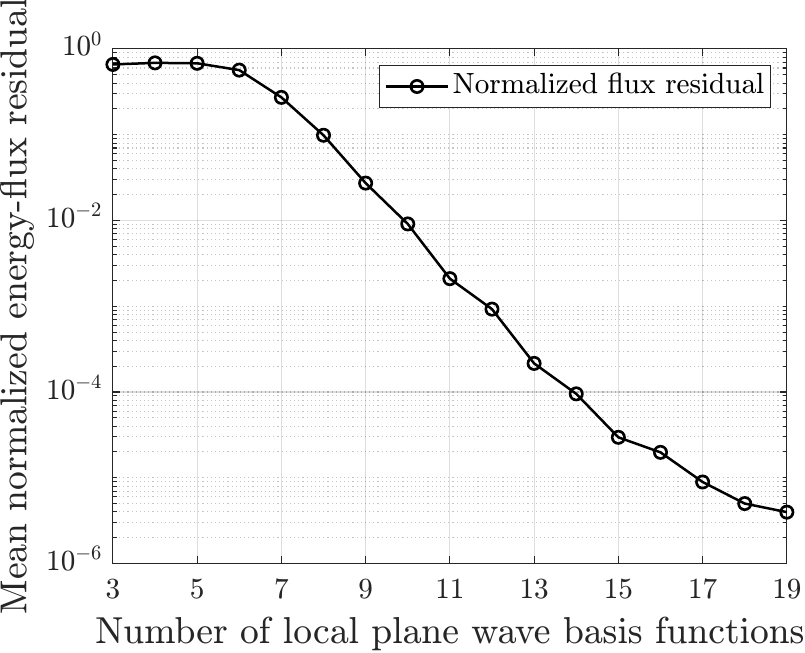}
\caption{Mean normalized energy-flux residual \eqref{eq:flux_check} of the T-matrix of the first layer of \S\ref{s:ex2}, 
against the number $p$ of local plane wave basis functions (semi-logarithmic scale).  The residual is evaluated for each selected propagating incident mode, normalized, and averaged in absolute value over these modes.}
\label{fig:flux_conv}
\end{figure}

To further assess the numerical consistency of the T-matrix approximation, we monitor the energy-flux balance in the first layer $\Omega_1$, which has real permittivity. 
For any incident vector $(\aup,\adown)$ and $(\bup,\bdown)=\mathbf \Tmat_1^M(\aup,\adown)$, the exact T-matrix satisfies the flux identity \eqref{eq:flux_check}.
We evaluate the left-hand side of \eqref{eq:flux_check} for each propagating incident mode, normalize it by $\sum_{n}\{ |\beta^+_{1,n}| (|b^\uparrow_{1,n}| ^2+|a^\downarrow_{1,n}|^2)+|\beta^-_{1,n} |( |b^\downarrow_{1,n}|^2 +|a^\uparrow_{1,n}|^2)\}$,
and average its absolute value over these modes. The truncation order is $M=15$. 
Figure \ref{fig:flux_conv} shows that this residual decreases as the number $p$ of local plane wave basis functions increases, 
which indicates that the computed T-matrix increasingly satisfies the expected energy balance.

\subsection{Comparison of the DtN-TDG with and without T-matrix}
\label{s:multi_layers}

We compare the DtN-TDG method on several layers, coupled with the T-matrix approach, against the same discretization applied treating the full domain $\Omega$ as a single layer, as in \cite{Monforte2025}.
We select a situation in which the stacked layers are all copies of a small number of layers: in this case the T-matrix method is even more convenient since the local T-matrix is computed just once per shape.

\begin{figure}[htbp]
\centering
\begin{tikzpicture}[every node/.style={scale=0.8}]
\draw (0,-1) -- (0,1) -- (4,1) -- (4,-1) -- (0,-1);
\draw (2,-0.5) -- (2.5,0.5) -- (3,-0.5) --  (2,-0.5);
\draw (0,1) node[left]{$1$};
\draw (0,-1) node[left]{$-1$};
\draw (0,-1) node[below]{$0$};
\draw (4,-1) node[below]{$4$};
\node[] at (1,0) {$\varepsilon = 1$};
\node[] at (2.5,0) {$\varepsilon_{in} = 2$};
 \end{tikzpicture}\hspace*{8pt}
 \begin{tikzpicture}[every node/.style={scale=0.8}]
 \draw (0,-1) -- (0,1) -- (4,1) -- (4,-1) -- (0,-1);
 \filldraw[gray!40!white] (1, 0) -- (1.5, 0.5) -- (2, 0) --  (1.5, -0.5) -- (1, 0);
 \draw (1, 0) -- (1.5, 0.5) -- (2, 0) --  (1.5, -0.5) -- (1, 0);
 \draw (0,1) node[left]{$1$};
 \draw (0,-1) node[left]{$-1$};
 \draw (0,-1) node[below]{$0$};
 \draw (4,-1) node[below]{$4$};
 \node[] at (3,0) {$\varepsilon = 1$};
 \end{tikzpicture}
 \caption{The reference layers for the multiple layer example in \S\ref{s:multi_layers}.}
 \label{fig:multiple_layers}
\end{figure}
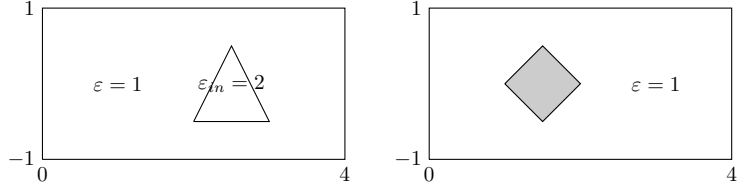

Let the domain $\Omega = (0,4) \times (-14, 14)$ be composed of $N=10$ layers, each of which is a vertical translation and possibly a reflection (recalling Remark~\ref{rem:Reflect}) of the two reference layers displayed in Figure \ref{fig:multiple_layers}.
These reference layers are the rectangle $(0,4)\times(-1,1)$ with $\varepsilon=1$ and either (i) a penetrable triangle with vertices $(2,-0.5)$, $(2.5,0.5)$ and $(3,-0.5)$, with $\varepsilon_{in}=2$ inside the triangle, and (ii) an impenetrable square with vertices $(1, 0)$, $(1.5, 0.5)$, $(2, 0)$ and $(1.5, -0.5)$.

We select an incident plane wave with wavenumber $k=8$ and incident angle $\theta=-\pi/4$.
We choose $h=0.5$ as the mesh parameter, $M=18$ as the Fourier modes truncation order, and select $p=20$ basis function in every element of the mesh.
Using the MATLAB ``timeit'' function, we compare the execution time of the two methods on multiple runs and compute the average: the execution time of the DtN-TDG method on the whole $\Omega$ is $22.5$ seconds on average, while the use of the T-matrix reduces the computational time to $4.2$ seconds.
The former scheme requires the solution of a DtN-TDG linear system of $16320 \times 16320$, the latter the solution of several systems of sizes $1640 \times 1640$ and $1440 \times 1440$ (for the two reference layers in Figure~\ref{fig:multiple_layers}).
Moreover, adding more layers does not considerably affect the computational cost of the T-matrix method, since most of the time is spent on approximating the local T-matrices. 
To compute the T-matrices in parallel we use the \texttt{parfor} MATLAB function, which uses the \texttt{Processes} cluster profile and $6$ parallel workers.

The execution time of the monolithic DtN-TDG method, $22.5$ seconds, is split into about $19$ seconds for the assembly of the linear system and about $4$ seconds for its solution with a direct solver. 
Reordering the unknowns reduces its half-bandwidth to $259$, so it has a narrow band structure, but this does not change the total time substantially, since the cost is dominated by the assembly; moreover, the MATLAB \texttt{backslash} solver already applies a reordering of this type. 
In the T-matrix method, the same assembly code is applied only to the two reference layers, so part of the gain is due to the reuse of repeated layers: the matrix associated to a layer is assembled and factorized once, however many times the layer appears.
Assembly and factorization can be performed in parallel, further reducing the time.

\begin{figure}[htbp]
\centering
\includegraphics[width=\textwidth]{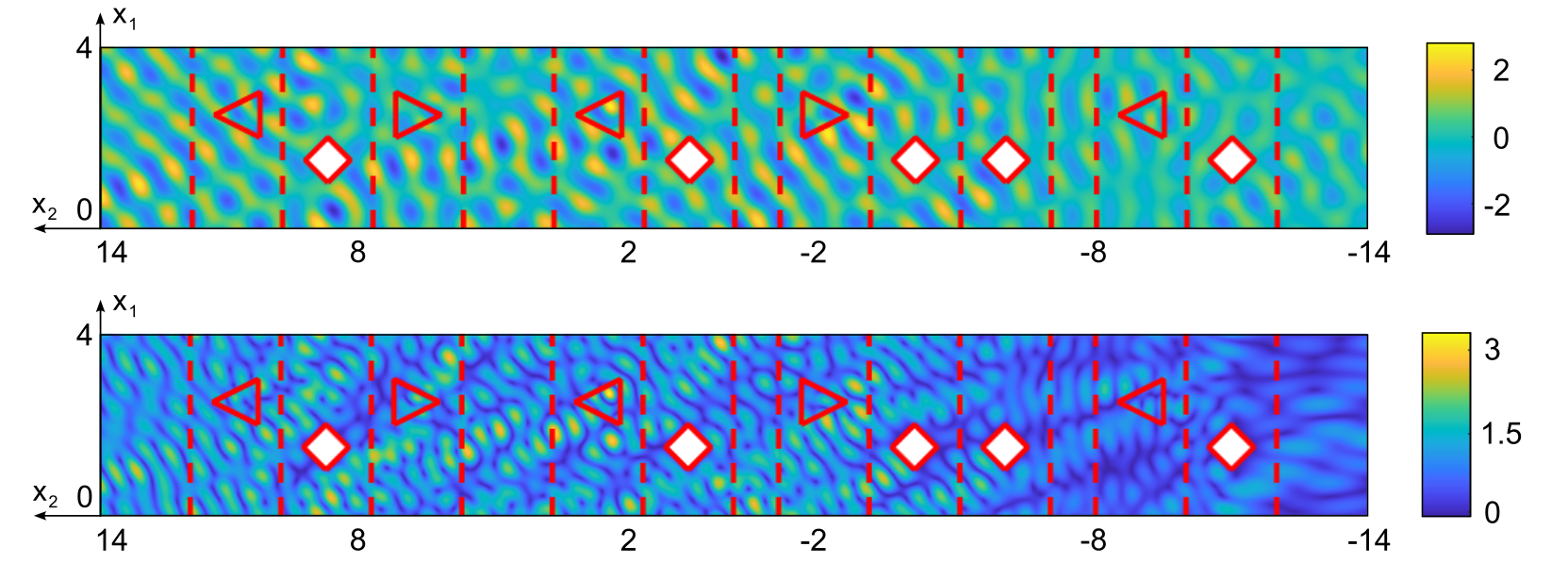}
\caption{The real part (top) and absolute value (bottom) of the T-matrix solution for the example in \S\ref{s:multi_layers}, with the axes rotated by $\pi/2$. Layers are delimited by dashed lines while the regions with different relative permittivity are separated by continuous lines.}
\label{fig:multi_layers}
\end{figure}

In Figure \ref{fig:multi_layers} we display the real part and absolute value of the solution obtained with the T-matrix approximation. The DtN-TDG method leads to the same result.

\subsection{Performance for increasing number \texorpdfstring{$N$}{N} of layers}
\label{s:ManyLayers}

We further investigate the scalability of the T-matrix approach by increasing the number of layers in the computational domain.
Starting from a configuration with $N=3$ layers, we iteratively add two layers at every step, one at the top and one at the bottom of the domain,
leaving no space between two adjacent layers, until reaching $N=45$ layers. 
If gaps between layers were present, then the T-matrix version of the method would be even more advantageous over the monolithic one.
Each new layer is selected randomly among the two reference configurations shown in Figure~\ref{fig:multiple_layers}, including their reflected versions described in Remark~\ref{rem:Reflect}.
The physical and numerical parameters are as in \S\ref{s:multi_layers}.
For every value of $N$, we measure the DtN-TDG computational time using the MATLAB commands \texttt{tic} and \texttt{toc}, for the following approaches: 
(i) the DtN-TDG method applied to the whole domain (without T-matrix); 
(ii) the T-matrix method, where the T-matrices of the two reference layers are computed once and reused for all layers, as in \S\ref{s:multi_layers}; 
(iii) the T-matrix method, where the repetition of the layers is not exploited and the T-matrix of each of the $N$ layers is computed, 
in serial;
(iv) same as (iii) but with the T-matrices computed in parallel (Matlab's \texttt{parfor}, 6 workers). 
Cases (iii) and (iv) have the same cost as for a structure made of $N$ different layers.

\begin{figure}[htbp]
\centering
\includegraphics[width=.7\textwidth]{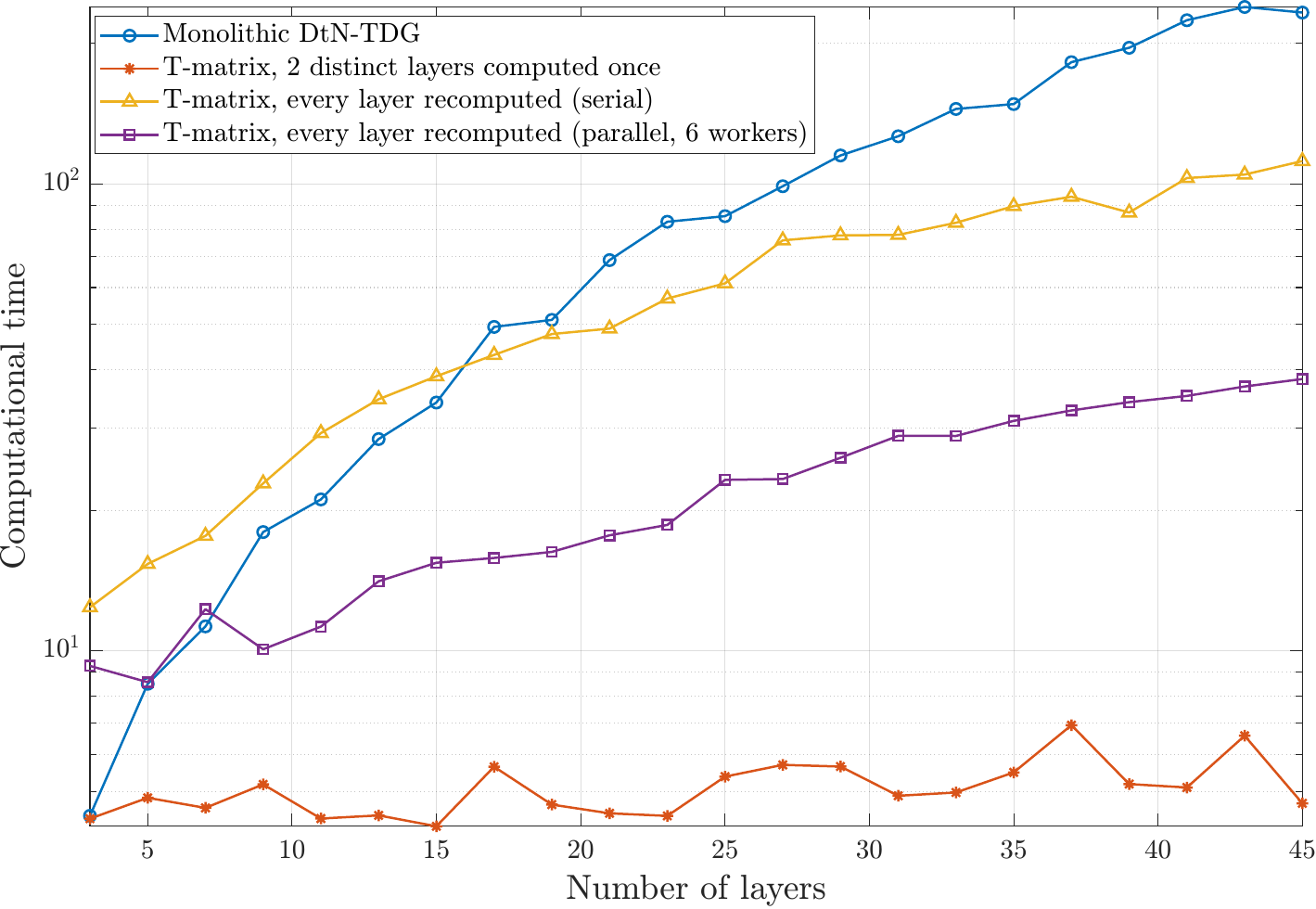}
\caption{Computational time for the example in \S\ref{s:ManyLayers}, plotted against the number of layers $N$, in semi-logarithmic scale, for the DtN-TDG method on the full domain and for the T-matrix accelerated approach: T-matrices of the two reference layers computed once and reused for all layers; T-matrix of every layer recomputed, in serial and in parallel (\texttt{parfor}, 6 workers).}
\label{fig:time_vs_layers}
\end{figure}

Figure~\ref{fig:time_vs_layers} displays the computational time as a function of the number of layers for both the standard ``monolithic'' DtN-TDG method applied to the full domain (blue circles) and the DtN-TDG combined with the T-matrix approach in the cases (ii) (red stars), (iii) (yellow triangles) and (iv) (purple squares) above.
As expected, the computational cost of the DtN-TDG discretization grows significantly with the number of layers, since the size of the global linear system increases proportionally to the domain height.
Indeed, at $N=45$ the mesh includes 2\,978 triangles and the DtN-TDG linear system has dimension 59\,560.
In contrast, the T-matrix method with the T-matrices computed once and reused exhibits a much milder growth in computational time, since most of the computational effort is associated with the approximation of the local T-matrices, which does not depend on the number of layers.
This curve measures the benefit of the presence of repeated layers and gives the lowest times. 
When the T-matrix of every layer is recomputed as if the $N$ layers were all distinct,
the time grows linearly with $N$, as expected from $N$ local problems of fixed size.
However, the computational time is still 
smaller than that needed for the monolithic TDG solver, since the $N$ TDG systems are of much smaller dimension when compared to the monolithic solver.
These results confirm that the T-matrix strategy is particularly advantageous in multilayer configurations, and even more so if the structures contained in the layers are repeated.

\subsection{Comparison with the RCWA method in proximity of BICs} \label{s:Comparison}

Among the numerical methods developed to simulate periodically patterned photonic structures, a well-established one is the Rigorous Coupled-Wave Analysis (RCWA), also known as Fourier Modal Method (FMM). This method, coupled with the scattering matrix propagation algorithm, is largely used to approximate the wave scattering on structures with periodicity in 1 or 2 dimensions, and layered on the third dimension \cite{Whittaker1999}. In this section, we use RCWA as a benchmark and compare its performance against the TDG T-matrix method. 

As a test structure, we consider a modified distributed Bragg reflector (DBR). An ordinary DBR would consist of an arrangement of alternating high-permittivity and low-permittivity layers, homogeneous in two dimensions, and would reflect waves whose wavelengths are contained in an interval known as bandgap \cite{Yariv}. However, if one replaces one or many homogeneous layers of the original DBR with patterned layer(s), additional photonic functionality promoted by unusual wave phenomena can be observed \cite{LI2026133360}. In particular, arrays of high-permittivity scatterers sandwiched between low-permittivity semi-infinite background regions support guided Rayleigh--Bloch waves above the light line, also known as bound states in the continuum (BIC). Such modes occur for isolated pairs in the domain $(k,\alpha_{0})$, where $k$ is the vacuum wavevector and $\alpha_0$ is the quasi-periodicity parameter, i.e.~the Bloch wavevector. Notably, a pure BIC is degenerate with a state belonging to a continuum of radiating states; nonetheless, the BIC cannot be excited from the far field as its eigenfunction is exponentially vanishing towards the far field. However, by slightly perturbing the wavevectors from the $m$-th BIC condition $(k,\alpha_{0}) = (k_m,\alpha_{0,m})$, the BIC becomes a \textit{quasi-BIC} and its presence is revealed as a sharp peak in the far-field response functions (reflectance, transmittance, diffraction efficiency, or absorptance). This can be understood as a resonance in the frequency domain (or equivalently in the $k$ domain), whose quality factor depends upon the Bloch wavevector; the quality factor decreases as $\alpha_{0}$ deviates from $\alpha_{0,m}$. In the literature, these resonances are known as leaky Rayleigh--Bloch waves, guided mode resonances, or quasi-normal modes, and play an increasingly relevant role in photonics \cite{Marinica2008,hsu_observation_2013,hsu_bound_2016,kang_applications_2023,Collin2014,lalanne_light_2018}.

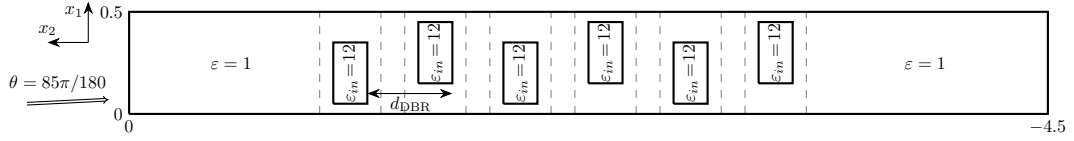
\begin{figure}[htbp]
	\centering
	\begin{tikzpicture}[scale=2.7,every node/.style={scale=0.67}]
		\draw[thick] (0,0) -- (0,0.5) -- (4.5,0.5) -- (4.5,0) -- (0,0);
		\draw[thick] (-1/12+13/12,0.05) -- (-1/12+13/12,0.35) -- (1/12+13/12,0.35) -- (1/12+13/12,0.05) -- (-1/12+13/12,0.05);
		\draw[thick] (-1/12+23/12,0.05) -- (-1/12+23/12,0.35) -- (1/12+23/12,0.35) -- (1/12+23/12,0.05) -- (-1/12+23/12,0.05);
		\draw[thick] (-1/12+11/4,0.05) -- (-1/12+11/4,0.35) -- (1/12+33/12,0.35) -- (1/12+33/12,0.05) -- (-1/12+33/12,0.05);
		\draw[thick] (-1/12+3/2,0.15) -- (-1/12+3/2,0.45) -- (1/12+3/2,0.45) -- (1/12+3/2,0.15) -- (-1/12+3/2,0.15);
		\draw[thick] (-1/12+7/3,0.15) -- (-1/12+7/3,0.45) -- (1/12+7/3,0.45) -- (1/12+7/3,0.15) -- (-1/12+7/3,0.15);
		\draw[thick] (-1/12+19/6,0.15) -- (-1/12+19/6,0.45) -- (1/12+19/6,0.45) -- (1/12+19/6,0.15) -- (-1/12+19/6,0.15);
		\draw[draw opacity=0.5,dashed] (-0.15+13/12,0) --  (-0.15+13/12,0.5);
		\draw[draw opacity=0.5,dashed] (0.15+13/12,0) --  (0.15+13/12,0.5);
		\draw[draw opacity=0.5,dashed] (-0.15+23/12,0) --  (-0.15+23/12,0.5);
		\draw[draw opacity=0.5,dashed] (0.15+23/12,0) --  (0.15+23/12,0.5);
		\draw[draw opacity=0.5,dashed] (-0.15+11/4,0) --  (-0.15+11/4,0.5);
		\draw[draw opacity=0.5,dashed] (0.15+11/4,0) --  (0.15+11/4,0.5);
		\draw[draw opacity=0.5,dashed] (-0.15+3/2,0) --  (-0.15+3/2,0.5);
		\draw[draw opacity=0.5,dashed] (0.15+3/2,0) --  (0.15+3/2,0.5);
		\draw[draw opacity=0.5,dashed] (-0.15+7/3,0) --  (-0.15+7/3,0.5);
		\draw[draw opacity=0.5,dashed] (0.15+7/3,0) --  (0.15+7/3,0.5);
		\draw[draw opacity=0.5,dashed] (-0.15+19/6,0) --  (-0.15+19/6,0.5);
		\draw[draw opacity=0.5,dashed] (0.15+19/6,0) --  (0.15+19/6,0.5);
		\draw (0,0) node[below]{$0$};
		\draw (4.5,0) node[below]{$-4.5$};
		\draw (0,0) node[left]{$0$};
		\draw (0,0.5) node[left]{$0.5$};
        \draw[{Stealth}-{Stealth}] (1/12+13/12,0.1) -- (1/12+3/2,0.1);
        \node[] at (1.38,0.05) {$d_{\text{DBR}}$};
        \node[rotate=90] at (13/12,0.2) {$\varepsilon_{in}\!=\!12$};
        \node[rotate=90] at (23/12,0.2) {$\varepsilon_{in}\!=\!12$};
        \node[rotate=90] at (33/12,0.2) {$\varepsilon_{in}\!=\!12$};
        \node[rotate=90] at (3/2,0.3) {$\varepsilon_{in}\!=\!12$};
        \node[rotate=90] at (7/3,0.3) {$\varepsilon_{in}\!=\!12$};
        \node[rotate=90] at (19/6,0.3) {$\varepsilon_{in}\!=\!12$};
        \node[] at (0.5,0.25) {$\varepsilon=1$};
        \node[] at (3.9,0.25) {$\varepsilon=1$};
        \node[] at (-0.35,0.15) {$\theta = 85\pi/180$};
        \draw[-{Implies},double] (-0.5 ,0.05) -- (-0.1,0.07);
        \draw[-{Stealth[slant=0]}] (-0.2 ,0.35) -- (-0.4,0.35);
        \draw[-{Stealth[slant=0]}] (-0.2 ,0.35) -- (-0.2,0.55);
        \node[] at (-0.27,0.5) {$x_1$}; 
        \node[] at (-0.4,0.41) {$x_2$}; 
	\end{tikzpicture}
	\caption{Sketch of the domain $\Omega$ for the example in \S\ref{s:Comparison}, with the axes rotated by $\pi/2$. The incident wave is represented by a double arrow.}
	\label{fig:BraggReflector}
\end{figure}

The structure under investigation is represented in Fig.~\ref{fig:BraggReflector}: here, the domain $\Omega = (0,0.5) \times (-4.5, 0)$ contains six identical high-permittivity ($\varepsilon_{in}=12$) rectangular obstacles of side lengths $0.3$ and $1/6$, placed along the $x_2$ direction with a DBR periodicity $d_{DBR}$.
The obstacles are alternately translated by $\pm 0.05$ along $x_1$ to lower the symmetry of the system and hence to promote the formation of quasi-BICs (q-BICs) \cite{koshelev_asymmetric_2018}. The system is periodic along $x_1$, and only a single period is represented in Fig.~\ref{fig:BraggReflector}. The exterior medium has relative permittivity $\varepsilon=1$. For the purpose of T-matrix TDG calculations, the domain is decomposed in layers, that are all translations of two reference layers $\Omega_1$ and $\Omega_2$. Each layer contains a single penetrable rectangle with relative permittivity $\varepsilon_{in}=12$. 

We select an incident plane wave with incident angle $\theta=-\pi\frac{85}{180}$ and vary the wavenumber $k$ (and consequently the wavelength) in the range $k = [0.3,1]\cdot2\pi$; in this range, only the zeroth diffraction order is allowed (i.e., the system operates in the metasurface regime).
To approximate the local T-matrices, we choose $h=0.04$ as the mesh parameter, $M=15$ as the Fourier modes truncation order, and select $p=15$ basis function in every element of the mesh, resulting in a mesh of $220$ triangles and a linear system of dimension $3300$ for every reference layer.
The quantity under analysis is the transmitted intensity, i.e.\ the square of the absolute value of the $0$-th Fourier coefficient $b_{6,0}^\downarrow$ of the solution $u$ on the lower face $\Gamma_6^-$ of the lowest layer.
The same system has been analyzed through a MATLAB-based RCWA implementation\footnote{PPML v4.0, available at \url{https://github.com/zan8simone/PPML_new}}, setting the Fourier truncation radius to 15 and the FFT discretization grid to 4096 points.

\begin{figure}[htbp]
	\centering
    \includegraphics[width=\textwidth]{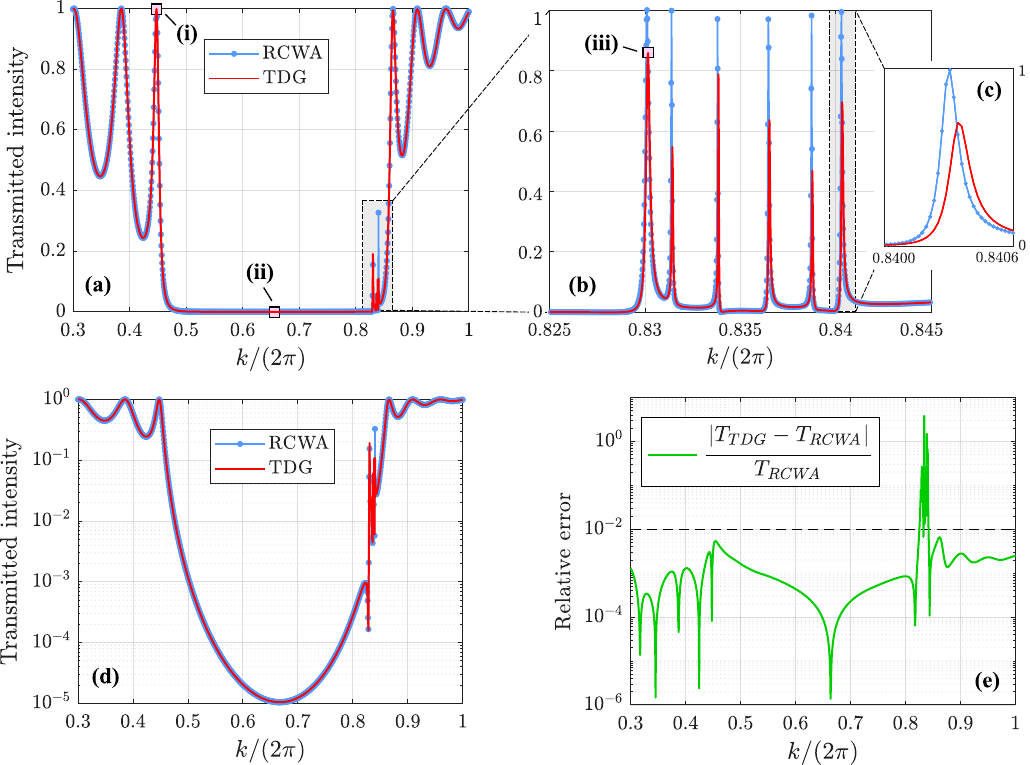}\hfill
	\caption{Panel (a-d): transmitted intensity $|b^{\downarrow}_{6,0}|^2$ of the grating of the example in Fig.~\ref{fig:BraggReflector} and \S\ref{s:Comparison}, plotted against $k/(2\pi)$ for the TDG T-matrix method and the RCWA method. Panels (b) and (c) are sequential zooms into the quasi-bound-states-in-continuum (qBIC) region. Note that in panel (a) the peaks in the gray region do not apparently reach unity because of undersampling. Panel (e): relative error of the transmitted intensity computed with the TDG method against the RCWA method. The error is always below $10^{-2}$ exception made for the qBIC region.}
	\label{fig:transmittted_intensity}
\end{figure}

The results are displayed in Figure~\ref{fig:transmittted_intensity}. The direct visual comparison in linear scale (panel a) and in logarithmic scale (panel d) show an excellent agreement, as both methods match in predicting the DBR bandgap extending in the interval $k=[0.5,0.8]\cdot 2\pi$, as well as the DBR sidelobes (i.e.~the transmittance oscillations outside the DBR bandgap).
As anticipated, the structure does not only support Bragg reflections, but also qBICs. Such resonances appear in the $k = [0.825,0.845]\cdot 2\pi$ region, as highlighted in panel (b).
Here it can be seen that both methods capture the very narrow resonances, except for minor shifts and peak height differences (panel c). Finally, panel (e) reports the relative error between the methods.
Except for the qBIC region, the relative error is always below $10^{-2}$ (dashed line), and often below $10^{-3}$. 

The discrepancy between TDG and RCWA can be ascribed to the fact that at exact BIC conditions the solution of the BVP \eqref{eq:BVP} is not unique, and the numerical method does not converge smoothly, as studied in \cite[\S5.2.1]{Monforte2025}. Additional studies on our system  revealed that the discrepancy between the two methods becomes larger as the incident angle approaches $\pi/2$, i.e.~as the qBICs approach the corresponding proper BICs. However, the case $\theta=85^{\circ}$ analyzed in Fig.~\ref{fig:BraggReflector} indicates that TDG is capable of properly capturing resonances with a quality factor of $10^4$, which is a significant result also in view of applications.

\begin{figure}[htbp]
	\centering
	\includegraphics[width=.99\textwidth]{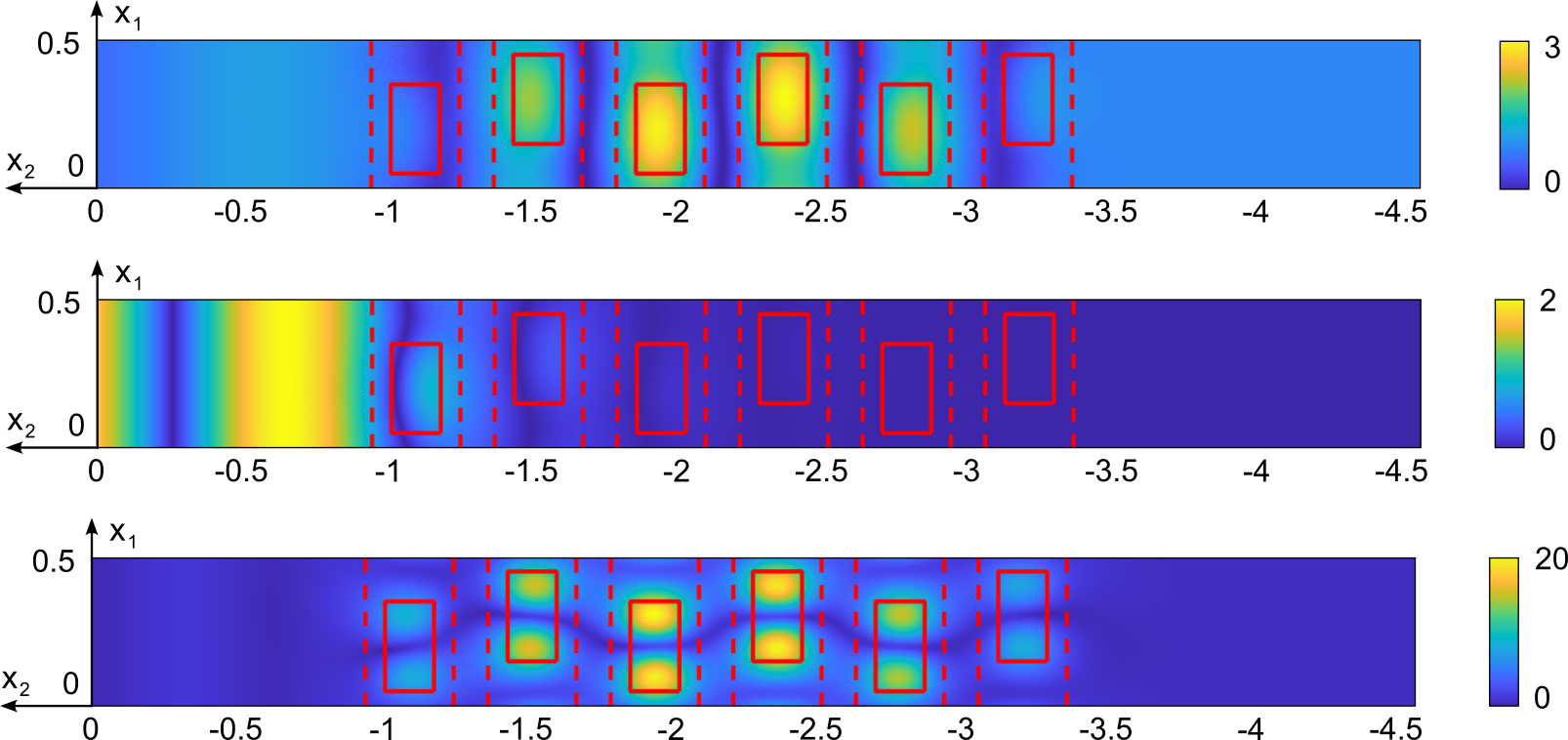}
	\caption{The absolute value of the T-matrix solution for the example in \S\ref{s:Comparison} in three different cases, with the axes rotated by $\pi/2$. Top: the wave weakly resonates inside the distributed Bragg reflector (DBR) and it is eventually transmitted [point labeled (i) in Fig.~\ref{fig:transmittted_intensity}].
    Middle: the wave is strongly reflected by the DBR  [point labeled (ii) in Fig.~\ref{fig:transmittted_intensity}].
    Bottom: the wave strongly resonates within the DBR thanks to guided mode resonance [point labeled (iii) in Fig.~\ref{fig:transmittted_intensity}].
    The layers are separated by dashed lines while the regions with different relative permittivity are separated by continuous lines.}
	\label{fig:DBR_fields}
\end{figure}

To complete the analysis of the modified DBR, we plot the numerical solution on $\Omega$ for three different wavenumbers $k$, corresponding to three different wave behaviors within the structure. The first case (Fig.~\ref{fig:DBR_fields}, top panel), corresponding to $k = 0.446446\cdot2\pi$, i.e.~point (i) in Fig.~\ref{fig:transmittted_intensity}, is the first low-frequency DBR sidelobe. Here, the wave is weakly resonating along the $x_2$ direction and eventually transmitted towards $x_2 \rightarrow -\infty$. The second case (Fig.~\ref{fig:DBR_fields}, middle panel), corresponding to $k = 0.65\cdot2\pi$, i.e.~point (ii) in Fig.~\ref{fig:transmittted_intensity}, is the DBR bandgap center; here, the wave is strongly reflected and the solution rapidly decreases within the structure. The third case (Fig.~\ref{fig:DBR_fields}, bottom panel), corresponding to $k = 0.8301\cdot2\pi$, i.e.~point (iii) in Fig.~\ref{fig:transmittted_intensity}, is the first qBIC. Here, the incident wave, traveling along $-x_2$, is scattered into waves traveling along $x_1$, and strongly resonates before being re-emitted towards the far field (note the color bar scale).
We have performed analogous calculations in RCWA, obtaining perfectly matching field profiles. 

This comparison shows that the T-matrix TDG is perfectly suitable to simulate structures of interest in applied photonics, yielding results comparable to those of a well-established method. 
Although the RCWA is faster than the TDG and appears to be slightly more accurate as it provides smoother spectra when exploring qBICs with quality factors larger than $\approx 10^4$, we highlight that TDG is much more general for what concerns the scatterers that can be treated.
Indeed, while RCWA only handles shapes with edges aligned to the coordinate axes, TDG can handle rotated inclusions (such as those in Fig.~\ref{fig:corner_sing} and \ref{fig:multiple_layers}) without extra cost; moreover, TDG supports impenetrable obstacles (Dirichlet boundary conditions), which is not feasible in RCWA.  

\section*{CRediT authorship contribution statement}
\textbf{Armando Maria Monforte}: Writing – Original Draft, Visualization, Validation, Software, Methodology, Investigation, Formal analysis, Conceptualization; \textbf{Andrea Moiola}: Writing – review \& editing, Supervision, Project administration, Methodology, Formal analysis, Conceptualization; \textbf{Simone Zanotto}: Writing – review \& editing, Software, Visualization, Validation, Investigation.

\section*{Acknowledgments}
Andrea Moiola acknowledges support by the PRIN project “ASTICE” (202292JW3F), funded by the European Union–NextGenerationEU. 

\section*{Declaration of competing interest}
The authors declare that they have no known competing financial interests or personal relationships that could have appeared to influence the work reported in this paper

\section*{Data availability}
The code and the replication data for all the experiments are available on \url{https://github.com/Arma99dillo/LayeredTMatrix}

\addcontentsline{toc}{section}{References}
\printbibliography

\end{document}